\documentclass[a4paper,11pt,fleqn]{article}
\usepackage{amsmath,amssymb,amsthm,graphicx,subfigure,float,caption,epstopdf,tabularx,color, bm,amsfonts,epic}
\usepackage[top=1in, bottom=1in, left=1.25in, right=1.25in]{geometry}
\usepackage{appendix}
\usepackage{multirow}
\usepackage{diagbox}
\usepackage{appendix}

\allowdisplaybreaks
\newtheorem{thm}{Theorem}[section]
\newtheorem{shn}{Scheme}[section]
\newtheorem{lem}{Lemma}[section]

\newtheorem{rmk}{Remark}[section]
\newtheorem*{prf}{Proof}
\numberwithin{equation}{section}

\usepackage{indentfirst}
\graphicspath{{Fig}}

\begin{document}
\title{Arbitrary high-order structure-preserving schemes for the generalized
Rosenau-type equation }
\author{Chaolong Jiang$^{1,2}$, Xu Qian$^1$\footnote{Correspondence author. Email:
qianxu@nudt.edu.cn.}, Songhe Song$^1$, and Chenxuan Zheng$^{2}$ \\
{\small $^1$ Department of Mathematics, College of Science, }\\
{\small National University of Defense Technology, Changsha, 410073, P.R. China}\\
{\small $^2$ School of Statistics and Mathematics, }\\
{\small Yunnan University of Finance and Economics, Kunming 650221, P.R. China}\\
}
\date{}
\maketitle

\begin{abstract}
In this paper, we are concerned with arbitrarily high-order momentum-preserving and energy-preserving schemes for solving the generalized
Rosenau-type equation, respectively. The derivation of the momentum-preserving schemes is made within the symplectic Runge-Kutta method, coupled with the standard Fourier pseudo-spectral method in space. Then, combined with the quadratic
auxiliary variable approach and the symplectic Runge-Kutta method, together with the standard Fourier pseudo-spectral method, we present a class of high-order mass- and energy-preserving schemes for the Rosenau equation. Finally, extensive numerical tests and comparisons are also addressed to illustrate the performance of the proposed schemes.  \\[2ex]
\end{abstract}

\section{Introduction}
In this paper, we consider the following generalized
Rosenau-type equation
\begin{align}\label{Ros-equation}
\left \{
 \aligned
&\partial_tu(x,t)+\kappa\partial_x u(x,t)-\delta\partial_{xxt}u(x,t)+ b \partial_{xxx}u(x,t)\\
&~~~~~~~~~~~~~~~~~~~~~+\alpha\partial_{xxxxt}u(x,t)+\beta \partial_x(u(x,t)^p)=0,\ {x}\in\Omega\subset\mathbb{R},\ t>0,\\
&u({x},0)=u_0({x}),\ {x}\in{\Omega}\subset\mathbb{R},
\endaligned
 \right.
  \end{align}
where $t$ is the time variable, ${x}$ is the spatial variable, $u:=u({x},t)$ is the real-valued wave function, $\kappa, \delta>0, b, \alpha>0$ and $\beta$ are given real constants, $p$ is a given positive integer, $u_0({x})$ is a given initial condition and $\Omega=[x_l,x_r]$ is a bounded domain. Here, the Rosenau equation \eqref{Ros-equation} will be augmented with periodic boundary condition.

 The Rosenau equation was originally introduced to describe the dynamics of dense discrete systems\cite{Rosenau1988}. Nowadays, it has played an important role in fluid mechanics as well as atmosphere and ocean. Moreover, when $u$ is assumed to be smooth, the equation \eqref{Ros-equation} satisfies the following Hamiltonian formulation
\begin{align}\label{Ros-Hamiltonian-system}
u_t=\mathcal J\frac{\delta\mathcal H}{\delta u},
\end{align}
where $\mathcal J = -(1-\delta\partial_{xx}+\alpha\partial_{xxxx})^{-1}\partial_x$ is a Hamiltonian operator and $\mathcal H$ is the Hamiltonian functional, namely
\begin{align}\label{Ros-Hamiltonian-energy}
\mathcal H(t) =\int_{\Omega}\Big(\frac{\kappa}{2}u^2-\frac{b}{2}u_x^2+\frac{\beta}{p+1} u^{p+1}\Big)d{x},\ t\ge0,\
\end{align}
Besides the Hamiltonian energy \eqref{Ros-Hamiltonian-energy}, the equation \eqref{Ros-equation} also conserves the mass 
\begin{align}\label{Ros-mass-energy}
\mathcal M(t)=\int_{\Omega}ud{x}\equiv \mathcal M(0),\ t\ge0,
\end{align}
 and the momentum
\begin{align}\label{Ros-momentum}
\mathcal I(t) =\int_{\Omega}\Big(\frac{1}{2}u^2+\frac{\delta}{2} u_x^2+\frac{\alpha}{2} u_{xx}^2\Big)d{x}\equiv \mathcal I(0),\ t\ge0.
\end{align}

In the numerical computation of such Hamiltonian partial differential equations, it is often preferable to design some special numerical schemes that inherit one or
more intrinsic properties of the original system exactly in a discrete sense; they are called structure-preserving schemes \cite{FQ10,FM2011,ELW06}. In \cite{ChungAA98},  Chung proposed an implicit finite difference (IFD) scheme, which can satisfy the discrete analogue of momentum \eqref{Ros-momentum}. Furthermore, they prove rigorously in mathematics that the scheme is
second-order accurate in time and space. Subsequently, Omrani et al.\cite{OAAKamc08} developed a linearly implicit momentum-preserving finite difference scheme for the classical Rosenau equation, in which a linear system is to be solved at every time step. Thus it is computationally much cheaper than that of the IFD scheme. Over the years, various momentum-preserving schemes for the equation \eqref{Ros-equation} have been proposed and analyzed (e.g., see Refs \cite{AQJMS2022,AOAA2015,HeND2016,LiCMA2016,PZAMM2012,WLWCMA2017,WZJCM2019,WDCAM2018,WDYAA2021,WCCPMMAS2021}). However, to the best of our knowledge, all of existing momentum-preserving schemes are only second-order accurate in time at most. It has been shown in \cite{GZW2020jcp,JCQSjsc2022} that, compared with the second-order schemes, the high-order ones not only provide smaller numerical errors as a large time step chosen, but also will be more advantages in the robustness. Consequently, the first aim of this paper is to present a novel paradigm for developing arbitrary high-order momentum-preserving schemes for the equation \eqref{Ros-equation}.

Besides the momentum conservation law \eqref{Ros-momentum}, the equation \eqref{Ros-equation} also satisfies the Hamiltonian energy \eqref{Ros-Hamiltonian-energy}, which is one of the most important first integrals of the Hamiltonian system. Based on the averaged vector field method \cite{QM08}, Cai et al. \cite{CLZCNSNS2018} proposed a second-order energy-preserving scheme, and two fourth-order energy-preserving schemes are then proposed by using the composition ideas \cite{ELW06}. Nevertheless, it is shown in \cite{ELW06} that  the high-order schemes obtained by the composition method will be at the price of a terrible zig-zag of the step points (see Fig. 4.2 in Ref. \cite{ELW06}), which may be tedious and time consuming. Thus, the construction of high-order energy-preserving schemes for the Rosenau equation \eqref{Ros-equation} seems to be still at its beginning stage. In this paper, the second aim is to present a new strategy for proposing arbitrary high-order energy-preserving schemes for the Rosenau equation \eqref{Ros-equation} based on the quadratic auxiliary variable (QAV) approach\cite{CGHW-NMTMA-2022,GHWW2022,Tapley-SISC2022}, which is inspired by the idea of the energy quadratization (EQ) approach \cite{SXY18,SXY19siamrev,YZW17}.

The rest of this paper is organized as follows. In Section \ref{Ros-section-2}, the high-order momentum-preserving and energy-preserving schemes for the equation \eqref{Ros-equation} are proposed, respectively and their structure-preserving properties are analysed in details. Extensive numerical examples and some comparisons are addressed to illustrate the performance of the proposed schemes in Section \ref{Ros-section-3}. We draw some conclusions in Section \ref{Ros-section-4}.


\section{High-order structure-preserving schemes}\label{Ros-section-2}
In this section, we will propose high-order momentum-preserving schemes and energy-preserving schemes for the equation \eqref{Ros-equation}, respectively.

To be the start, let the spatial step size $h=\frac{x_r-x_l}{N}$ with an even positive integer $N$, and denote the grid points by $x_{j}=jh$ for $j=0,1,2,\cdots,N$; set $u_{j}$ be the numerical approximation of  $u(x_j,t)$ for $j=0,1,\cdots,N$, and $U:=(u_{0},u_{1},\cdots,u_{N-1})^T$ be the solution vector; we also define the discrete inner product, $l^2$-norm and $l^{\infty}$-norm as, respectively,
\begin{align*}
&\langle U,V\rangle_{h}=h\sum_{j=0}^{N-1}u_{j}v_{j},\ \|U\|_{h}^2=\langle U,U\rangle_{h},\ \|{U}\|_{h,\infty}=\max\limits_{0\le j<N-1}|u_{j}|.
\end{align*}
  In addition, we denote $`\cdot$' as the element product of vectors $U$ and  $V$, that is
\begin{align*}
U\cdot V=&\big(u_{0}v_{0},u_{1}v_{1},\cdots,u_{N-1}v_{N-1}\big)^{T}.
\end{align*}
For brevity, we denote ${U}\cdot U$ as $U^2$.

As achieving high order accuracy in time, the spatial accuracy shall be comparable to that of the time-discrete discretization. Actually, we consider the periodic boundary condition in this paper, so that the Fourier pseudo-spectral method is a very good choice
because of the high order accuracy and the fast Fourier transform (FFT) algorithm. Thus, we first expounded the Fourier pseudo-spectral method, as follows.

Let the interpolation space as
\begin{align*}
&\mathcal S_h=\text{span}\{l_{j}(x),\ 0\leq j\leq N-1\}
\end{align*}
where $l_{j}(x)$ is trigonometric polynomials of degree $N/2$ given by
\begin{align*}
  l_{j}(x)=\frac{1}{N}\sum_{k=-N/2}^{N/2}\frac{1}{c_{k}}e^{\text{i}k\mu(x-x_{j})},\ c_{k}=\left \{
 \aligned
 &1,\ |k|<\frac{N}{2},\\
 &2,\ |k|=\frac{N}{2},
 \endaligned
 \right.
\end{align*}
with $\mu=\frac{2\pi}{x_r-x_l}$. We then define the interpolation operator $\mathcal{I}_{N}: C(\Omega)\to \mathcal S_h$ as \cite{CQ01,CW-2016-Siam}
\begin{align*}
\mathcal{I}_{N}u(x,t)=\sum_{k=0}^{N-1}u_{k}(t)l_{k}(x),
\end{align*}
where $u_{k}(t)=u(x_{k},t),\ k=0,1,2,\cdots,N-1$. 

Taking the partial derivative with respect to $x$ at the collocation points $x_{j}$, we have
\begin{align}\label{Ros-fourier-2.2}
\frac{\partial^{r} \mathcal{I}_{N}u(x_{j},t)}{\partial x^{r}}
&=\sum_{k=0}^{N-1}u_{k}(t)\frac{d^{r}l_{k}(x_{j})}{dx^{r}}=\sum_{k=0}^{N-1}(D_r)_{j,k}u_{k}(t),\ j=0,1,\cdots, N-1,
\end{align}
where $D_r$ represents the spectral differential matrix with elements given by
\begin{align*}
(D_r)_{j,k}= \frac{d^{r}l_{k}(x_{j})}{dx^{r}},\ j,k=0,1,\cdots,N-1.
\end{align*}
In particular, we have\cite{JCWLMaxwell-2017,ST06}
\begin{align*}
&(D_{1})_{j,k}=
 \left \{
 \aligned
 &\frac{1}{2}\mu (-1)^{j+k} \cot(\mu \frac{x_{j}-x_{k}}{2}),\ &j\neq k,\\
 &0,\quad \quad \quad \quad \quad \quad \quad \ \ \ ~&j=k,
 \endaligned
 \right.
\end{align*}
\begin{align*}
 &(D_{2})_{j,k}=
 \left \{
 \aligned
 &\frac{1}{2}\mu^{2} (-1)^{j+k+1}\csc^{2}(\mu \frac{x_{j}-x_{k}}{2}),\ &j\neq k,\\
 &-\mu^{2}\frac{N^{2}+2}{12},\quad \quad \quad \quad ~ &j=k,
 \endaligned
 \right.
 \end{align*}
\begin{align*}
 &(D_{3})_{j,k}=
 \left \{
 \aligned
 &\frac{3\mu^{3}}{4} (-1)^{j+k}\cos(\mu \frac{x_{j}-x_{k}}{2})\csc^{3}(\mu \frac{x_{j}-x_{k}}{2})\\
 &~~~~~~~~~~~~~~~~~~~+\frac{\mu^{3}N^{2}}{8}(-1)^{j+k+1}\cot(\mu \frac{x_{j}-x_{k}}{2}),\ &j\neq k,\\
 &0,\quad \quad \quad \quad \quad \quad \quad~ ~~~~~~~~~~~~~~~~~~~~~~~~~~&j=k,
 \endaligned
 \right.
 \end{align*}
and
\begin{align*}
 &(D_{4})_{j,k}=
 \left \{
 \aligned
 &\mu^{4}(-1)^{j+k}\csc^{2}(\mu \frac{x_{j}-x_{k}}{2})\Big(\frac{N^2}{4}-\frac{1}{2}-\frac{3}{2}\cot^{2}(\mu \frac{x_{j}-x_{k}}{2})\Big),\ &j\neq k,\\
 &\mu^{4}(\frac{N^{4}}{80}+\frac{N^{2}}{12}-\frac{1}{30}),\quad \quad \quad \quad \quad ~~~~~~~~~~~~~ &j=k.
 \endaligned
 \right.
 \end{align*}

   \begin{rmk}\label{Ros-remark-2.1} We should note that \cite{GCW14,ST06}
\begin{align*}
& {D}_{r}=
 \left \{
 \aligned
 &\mathcal{F}_{N}^{H}\Lambda^{r}\mathcal{F}_{N},\ \ \ r\ \text{is an odd integer,}\\
 &\mathcal{F}_{N}^{H}\tilde{\Lambda}^{r}\mathcal{F}_{N},\ \ \ r\ \text{is an even integer,}
 \endaligned
 \right.
 \end{align*}
where $\Lambda$ and $\tilde{\Lambda}$ are
\begin{align*}
&\Lambda=\text{\rm diag}\Big(\text{\rm i}\mu[0,1,\cdots,\frac{N}{2}-1,0,-\frac{N}{2}+1,\cdots,-2,-1]\Big),\\
&\tilde{\Lambda}=\text{\rm diag}\Big(\text{\rm i}\mu[0,1,\cdots,\frac{N}{2}-1,\frac{N}{2},-\frac{N}{2}+1,\cdots,-2,-1]\Big),
\end{align*}
and $\mathcal F_N$ is the discrete Fourier transform (DFT) and $\mathcal F_N^H$ represents the conjugate
transpose of $\mathcal F_N$.
   \end{rmk}

We set $t_{n}=n\tau,$ and $t_{ni}=t_n+c_i\tau,\ i=1,2,\cdots,s$, $n = 0,1,2,\cdots,$ where $\tau$ is the time step size and $c_1,c_2,\cdots,c_s$ are distinct real numerbers (usually $0\le c_i\le 1$). The approximations of the function $u(x,t)$ at points  $(x_j,t_n)$ and $(x_j,t_{ni})$ are denoted by $u_j^n$ and $u_{j}^{ni}$, respectively.

\subsection{High-order momentum-preserving scheme}\label{Ros-section-2.1}
In this section, we propose a class of high-order momentum-preserving schemes for the equation \eqref{Ros-equation}. To be the start, we rewrite \eqref{Ros-equation} into
\begin{align}\label{Ros-momentum-equation}
\mathcal Au_t=\mathcal F(u)u,\ \mathcal F(u)=-\big[\kappa\partial_x+b\partial_{xxx}
+\frac{p\beta}{p+1}(u^{p-1}\partial_x+\partial_x(u^{p-1}))\big],
\end{align}
where $\mathcal A =1-\delta\partial_{xx}+\alpha\partial_{xxxx}$ is a self-adjoint operator, $\mathcal F(u)$ is an anti-adjoint operator and $\mathcal F(u)u$ is defined by
 \begin{align*}
 \mathcal F(u)u=-\big[\kappa\partial_xu+b\partial_{xxx}u
+\frac{p\beta}{p+1}(u^{p-1}\partial_xu+\partial_x(u^{p}))\big].
 \end{align*}
  The Fourier pseudo-spectral method is then employed to solve \eqref{Ros-momentum-equation} in space and we obtain
   \begin{align}\label{semi-Ros-momentum-equation}
\left \{
 \aligned
 &\mathcal A_h\frac{d}{dt}U=\mathcal F_h(U)U,\\
 & \mathcal F_h(U)=-\big[\kappa D_1+b D_3+\frac{p\beta}{p+1}\big(\text{diag}(U^{p-1})D_1+D_1 \text{diag}(U^{p-1})\big)\big],
 \endaligned
 \right.
\end{align}
where $\mathcal A_h =I-\delta D_2+\alpha D_4$ is a symmetric matrix, $\mathcal F_h(U)$ is anti-symmetric for $U$, and $\mathcal F_h(U)U$ is defined by
 \begin{align}\label{Ros-F-h-u}
 \mathcal F_h(U)U=-\big[\kappa D_1 U+bD_3U
+\frac{p\beta}{p+1}(U^{p-1}\cdot D_1U+D_1(U^{p}))\big].
 \end{align}
   \begin{thm} The semi-discrete system \eqref{semi-Ros-momentum-equation} preserves the following semi-discrete momentum conservation law
  \begin{align}\label{sm-SAV-energy-conservetion-law}
  \frac{d}{dt} I_h(t)=0,\ I_h(t)=\frac{1}{2}\langle\mathcal A_hU,U\rangle_h,\ t\ge0.
  \end{align}
   \end{thm}
   \begin{prf} With noting the symmetric property of $\mathcal A_h$ and anti-symmetric property of $\mathcal F_h(U)$ for $U$, we have
   \begin{align*}
\frac{d}{dt} I_h(t)= \langle\mathcal A_h\frac{d}{dt}U,U\rangle_h=\langle\mathcal F_h(U)U,U\rangle_h=0,
\end{align*}
which completes the proof.\qed
   \end{prf}
\begin{thm} If $p=2$, the semi-discrete system \eqref{semi-Ros-momentum-equation} preserves the following semi-discrete mass
  \begin{align}\label{sm-SAV-energy-conservetion-law}
  \frac{d}{dt} M_h(t)=0,\ M_h(t)=\langle U,1\rangle_h,\ t\ge0.
  \end{align}
   \end{thm}
 \begin{prf}It follows from \eqref{semi-Ros-momentum-equation} that
 \begin{align}
 \frac{d}{dt} M_h(t)=\langle\mathcal A_h^{-1}\mathcal F_h(U)U,{\bf 1}\rangle_h=\langle \mathcal F_h(U)U,\mathcal A_h^{-1} {\bf 1}\rangle_h=\langle \mathcal F_h(U)U,{\bf 1}\rangle_h.
 \end{align}
 With \eqref{Ros-F-h-u}, we have
 \begin{align}
 \langle \mathcal F_h(U)U,{\bf 1}\rangle_h &=  -\langle\kappa D_1 U+bD_3U
+\frac{p\beta}{p+1}(U^{p-1}\cdot D_1U+D_1(U^{p})),{\bf 1}\rangle_h\nonumber\\
&=-\frac{p\beta}{p+1}\langle D_1U,U^{p-1}\rangle_h.
 \end{align}
 As $p=2$, we obtain from the above equation
 \begin{align*}
 \langle \mathcal F_h(U)U,{\bf 1}\rangle_h=0,
 \end{align*}
 which implies that
 \begin{align*}
 \frac{d}{dt} M_h(t)=0.
 \end{align*}
 This completes the proof.\qed
 \end{prf}
 \begin{rmk} If $p>2$ and is a positive integer, we can deduce that the Fourier spectral differential matrix $D_1$ cannot satisfy the discrete equation
 \begin{align}\label{LR-2.8}
 \langle D_1 U,U^{p-1}\rangle_h=0, \ \text{\rm for $\forall$\ $U$.}
\end{align}
 Thus, the system \eqref{semi-Ros-momentum-equation} cannot conserve the semi-discrete mass \eqref{sm-SAV-energy-conservetion-law}, as $p>2$.
 \end{rmk}

We then apply an RK method to the system \eqref{semi-Ros-momentum-equation} in time to give a class of fully discrete schemes for the equation \eqref{Ros-equation}, as follows:
 \begin{shn}\label{Ros-scheme-2.1} Let $b_i,a_{ij}(i,j=1,\cdots,s)$ be real numbers and let $c_i=\sum_{j=1}^sa_{ij}$. For the given $U^n$, an s-stage Runge-Kutta method is given by
\begin{align}\label{Ros-momentum-schemes}
\left\{
\begin{aligned}
&\mathcal A_hK_i^n=\mathcal F_h(U^{ni})U^{ni},\ U^{ni}=U^n+\tau\sum_{j=1}^{s}a_{ij}K_j^n,\ i=1,2,\cdots,s,\\
&U^{n+1}=U^n+\tau\sum_{i=1}^{s}b_iK_i^n.
\end{aligned}
\right.
\end{align}
\end{shn}

\begin{thm} \label{Ros-equation-th2.1}  If the coefficients of the RK method satisfy
\begin{align}\label{Ros-RK-coeff}
b_ia_{ij}+b_ja_{ji}=b_ib_j, \ \forall \ i,j=1,\cdots,s,
\end{align}
 {\bf Scheme 2.1} conserves the following discrete momentum
\begin{align}\label{Ros-equation-new-2.5}
 I_h^{n+1}=I_h^n, \ I_h^n=\frac{1}{2}\langle U^n,\mathcal A_hU^n\rangle_h,\  n=0,1,2,\cdots.
\end{align}
\end{thm}
\begin{prf} It follows from the second equality of \eqref{Ros-momentum-schemes} that
\begin{align}\label{MP-Ros-equation-2.4}
 &I_h^{n+1}-{I}_h^n\nonumber\\
 &=\frac{1}{2}\langle U^{n+1},\mathcal A_hU^{n+1}\rangle_h-\frac{1}{2}\langle U^{n},\mathcal A_hU^{n}\rangle_h\nonumber\\
&=\frac{1}{2}\langle U^n+\tau\sum_{i=1}^{s}b_iK_i^n,\mathcal A_h(U^n+\tau\sum_{j=1}^{s}b_jK_j^n)\rangle_h-\frac{1}{2}\langle U^{n},\mathcal A_hU^{n}\rangle_h\nonumber\\
&=\frac{\tau}{2}\sum_{i=1}^{s}b_i\langle U^n,\mathcal A_hK_i^n\rangle_h+\frac{\tau}{2}\sum_{i=1}^{s}b_i\langle K_i^n,\mathcal A_hU^n\rangle_h+\frac{\tau^2}{2}\sum_{i,j=1}^{s}b_ib_j\langle K_i^n,\mathcal A_hK_j^n\rangle_h.
  \end{align}
  With noting
  \begin{align}\label{MP-Ros-equation-2.5}
  \frac{\tau}{2}\sum_{i=1}^{s}b_i\langle U^n,\mathcal A_hK_i^n\rangle_h&=\frac{\tau}{2}\sum_{i=1}^{s}b_i\langle U^{ni}-\tau\sum_{j=1}^{s}a_{ij}K_j^n,\mathcal A_hK_i^n\rangle_h\nonumber\\
  &=\frac{\tau}{2}\sum_{i=1}^{s}b_i\langle U^{ni},\mathcal A_hK_i^n\rangle_h-\frac{\tau^2}{2}\sum_{i,j=1}^{s}b_ja_{ji}\langle K_i^n,\mathcal A_hK_j^n\rangle_h.
  \end{align}
  Similarly, we have
  \begin{align}\label{MP-Ros-equation-2.6}
  \frac{\tau}{2}\sum_{i=1}^{s}b_i\langle K_i^n,\mathcal A_hU^n\rangle_h=\frac{\tau}{2}\sum_{i=1}^{s}b_i\langle K_i^n,\mathcal A_h U^{ni}\rangle_h-\frac{\tau^2}{2}\sum_{i,j=1}^{s}b_ia_{ij}\langle K_i^n,\mathcal A_hK_j^n\rangle_h.
  \end{align}
  We insert \eqref{MP-Ros-equation-2.5} and \eqref{MP-Ros-equation-2.6} into \eqref{MP-Ros-equation-2.4} and then use the symmetry of $\mathcal A_h$ to obtain
  \begin{align*}
 I_h^{n+1}-{I}_h^n&=\frac{\tau}{2}\sum_{i=1}^{s}b_i\big[\langle U^{ni},\mathcal A_hK_i^n\rangle_h+\langle K_i^n,\mathcal A_hU^{ni}\rangle_h\big]\\
 &~~~~~~~~~~~~~~~~~~~+\frac{\tau^2}{2}\sum_{i,j=1}^{s}( b_ib_j-b_ia_{ij}-b_ja_{ji})\langle K_i^n,\mathcal A_hK_j^n\rangle_h\\
  &=\tau\sum_{i=1}^{s}b_i\langle U^{ni},\mathcal A_hK_i^n\rangle_h +\frac{\tau^2}{2}\sum_{i,j=1}^{s}(b_ib_j-b_ia_{ij}-b_ja_{ji})\langle K_i^n,\mathcal A_hK_j^n\rangle_h.
  \end{align*}
  {The condition \eqref{Ros-RK-coeff} together with the equality}
\begin{align*}
\langle U^{ni},\mathcal A_hK_i^n\rangle_h=\langle U^{ni},\mathcal F_h(U^{ni})U^{ni}\rangle_h=0
\end{align*}
 implies $\ I_h^{n+1}={I}_h^n$. This completes the proof. \qed
\end{prf}

\begin{thm}\label{Ros-equation-th2.2} As $p=2$, {\bf Scheme 2.1} conserves the discrete mass, as follows:
\begin{align}\label{Ros-equation-new-2.16}
 M_h^{n+1}=M_h^n, \ M_h^n=\langle U^n,{\bf 1}\rangle_h,\  n=0,1,2,\cdots.
\end{align}
\end{thm}
\begin{prf} It follows from \eqref{Ros-momentum-schemes} that
\begin{align}
M_h^{n+1}&-M_h^n\nonumber\\
&=\tau\sum_{i=1}^{s}b_i\langle K_i^n,{\bf 1}\rangle_h\nonumber\\
&=\tau\sum_{i=1}^{s}b_i\langle \mathcal A_h^{-1}\mathcal F_h(U^{ni})U^{ni},{\bf 1}\rangle_h\nonumber\\
&=\tau\sum_{i=1}^{s}b_i\langle\mathcal F_h(U^{ni})U^{ni},{\bf 1}\rangle_h\nonumber\\
&=-\tau\sum_{i=1}^{s}b_i\langle\kappa D_1 U^{ni}+bD_3U^{ni}
+\frac{p\beta}{p+1}((U^{ni})^{p-1}\cdot D_1U^{ni}+D_1((U^{ni})^{p})),{\bf 1}\rangle_h\nonumber\\
&=-\frac{p\beta}{p+1}\langle D_1U^{ni},(U^{ni})^{p-1}\rangle_h.\nonumber
\end{align}
With noting $p=2$ and the antisymmetry of $D_1$, we have
\begin{align*}
\langle D_1U^{ni},(U^{ni})^{p-1}\rangle_h=0,
\end{align*}
which implies \eqref{Ros-equation-new-2.16}. This completes the proof. \qed
\end{prf}

\begin{rmk}\label{Ros-rmk-2.1} Assume that the initial condition $u_0({x})$ is sufficiently smooth, then it follows from \eqref{Ros-equation-new-2.5} that the numerical solution of {\bf Scheme 2.1} satisfies
\begin{align*}
&\sqrt{\|U^n\|_h^2+\delta\langle -D_2U^n,U^n\rangle_h+\alpha\langle D_4U^n,U^n\rangle_h}\\
&~~~~~~~~~~~~~~~~~~~~~~~=\sqrt{\|U^0\|_h^2+\delta\langle -D_2U^0,U^0\rangle_h+\alpha\langle D_4U^0,U^0\rangle_h}\le C,
\end{align*}
which implies that (noting $\alpha>0$ and $\delta>0$)
\begin{align*}
\|U^n\|_h\le C,\ \langle -D_2U^n,U^n\rangle_h\le C,\ \langle D_4U^n,U^n\rangle_h\le C,
\end{align*}
is uniformly bounded. Thus, {\bf Scheme 2.1} is unconditionally stable.
\end{rmk}

\begin{rmk}\label{Ros-RMK-2.4} {If we take $c_1,c_2,\cdots,c_s$ as the zeros of the $s$th shifted Legendre polynomial
\begin{align*}
\frac{d^s}{dx^s}\Big(x^s(x-1)^s\Big),
\end{align*}
the RK (or collocation) method based on these nodes has the order $2s$ and satisfies  the condition \eqref{Ros-RK-coeff} (see Refs. \cite{ELW06,Sanzs88,SCbook94} and references therein). In particular, the RK coefficients for $s=2$ and $s=3$ (denoted by 4th-Gauss method and 6th-Gauss method, respectively) are  given in Table \ref{Gauss-cllocation-method}} (e.g., see Ref. \cite{ELW06}). In addition, we introduce two notations, as follows:
\begin{itemize} \item 4th-MPS: using the 4th-Gauss method to {\bf Scheme 2.1};
\item 6th-MPS: using the 6th-Gauss method to {\bf Scheme 2.1}.
\end{itemize}
\begin{table}[H]
\centering
\begin{tabular}{c|cc}
${c}$ & ${A}$  \\
\hline
& ${b}^{T}$ \\
\end{tabular}
=
\begin{tabular}{c|cc}
$\frac{1}{2}-\frac{\sqrt{3}}{6}$ &$\frac{1}{4}$ & $\frac{1}{4}- \frac{\sqrt{3}}{6}$\\
$\frac{1}{2}+\frac{\sqrt{3}}{6}$ &$\frac{1}{4}+ \frac{\sqrt{3}}{6}$ &$\frac{1}{4}$ \\
\hline
                                 &$\frac{1}{2}$&$\frac{1}{2}$
\end{tabular}\\
\vspace{1mm}
\begin{tabular}{c|cc}
${c}$ & ${A}$  \\
\hline
& ${b}^{T}$ \\
\end{tabular}
=\begin{tabular}{c|ccc}
$\frac{1}{2}-\frac{\sqrt{15}}{10}$ &$\frac{5}{36}$ &  $\frac{2}{9}-\frac{\sqrt{15}}{15}$    &$\frac{5}{36}- \frac{\sqrt{15}}{30}$\\
$\frac{1}{2}$   &$\frac{5}{36}+ \frac{\sqrt{15}}{24}$  &  $\frac{2}{9}$  & $\frac{5}{36}- \frac{\sqrt{15}}{24}$ \\
$\frac{1}{2}+\frac{\sqrt{15}}{10}$ & $\frac{5}{36}+ \frac{\sqrt{15}}{30}$ & $\frac{2}{9}+\frac{\sqrt{15}}{15}$  & $\frac{5}{36}$ \\
\hline
& $\frac{5}{18}$ & $\frac{4}{9}$  & $\frac{5}{18}$
\end{tabular}
\caption{Gauss methods of order 4 and 6.}\label{Gauss-cllocation-method}
\end{table}
\end{rmk}


\subsection{High-order energy-preserving scheme}\label{Ros-section-2.2}

In this section, we propose a class of high-order energy-preserving schemes for the equation \eqref{Ros-equation}. Inspired by \cite{GHWW2022,JWG19,Tapley-SISC2022}, we first shall introduce appropriate quadratic auxiliary variables to reformulate the Hamiltonian energy into
a quadratic form. For clarity, we take $p=2,3$ and $5$ as examples to expound this procedure, as follows:
\begin{itemize}
\item {\bf Case I}: when $p=2$, we set
\begin{align}\label{2-IEQ-RLW-energy}
q:=q(x,t)=u^2.
\end{align}
Then, the Hamiltonian energy \eqref{Ros-Hamiltonian-energy} is rewritten into
\begin{align}
\mathcal H(t) =\int_{\Omega}\Big(\frac{\kappa}{2}u^2-\frac{b}{2}u_x^2+\frac{\beta}{3} uq\Big)d{x},\ t\ge0,
\end{align}
and according to the energy variational principle, we obtain the following reformulated system from \eqref{Ros-Hamiltonian-system}
\begin{align}\label{2-Ros-IEQ-reformulation}
\left\{
\begin{aligned}
&\partial_t u=\mathcal J\bigg (\kappa u+b u_{xx}+\frac{\beta}{3}q+ \frac{2\beta}{3}u^2\bigg), \\
&\partial_t q=2u\cdot\partial_tu,
\end{aligned}
\right.
\end{align}
 with the consistent initial conditions
\begin{align}\label{intial-Ros-IEQ-reformulation-2.18}
u({x},0)=u_0({x}),\ q({x},0)=(u_0({x}))^2.
\end{align}
\item {\bf Case II}: for $p=3$, the quadratic auxiliary variable is introduced as \eqref{2-IEQ-RLW-energy}, and the quadratic energy is given by
\begin{align}\label{3-IEQ-RLW-energy}
\mathcal H(t) =\int_{\Omega}\Big(\frac{\kappa}{2}u^2-\frac{b}{2}u_x^2+\frac{\beta}{4} q^2\Big)d{x},\ t\ge0,
\end{align}
which implies that the reformulated system is given by
\begin{align}\label{Ros-IEQ-reformulation}
\left\{
\begin{aligned}
&\partial_t u=\mathcal J\bigg (\kappa u+b u_{xx}+\beta uq\bigg), \\
&\partial_t q=2u\cdot\partial_tu,
\end{aligned}
\right.
\end{align}
 with the consistent initial conditions \eqref{intial-Ros-IEQ-reformulation-2.18}.

\item {\bf Case III}: As $p=5$, we start with introducing auxiliary variables
\begin{align}
q_1:=q_1(x,t)=u^2, q_2:=q_2(x,t)=uq_1,
\end{align}
 and then rewrite the Hamiltonian energy \eqref{Ros-Hamiltonian-energy} as
\begin{align}\label{IEQ-RLW-energy1}
\mathcal H(t) =\int_{\Omega}\Big(\frac{\kappa}{2}u^2-\frac{b}{2}u_x^2+\frac{\beta}{6} q_2^2\Big)d{x},\ t\ge0.
\end{align}
Similarly, we obtain reformulated system from \eqref{Ros-Hamiltonian-system}, as follows:
\begin{align}\label{Ros-IEQ-reformulation1}
\left\{
\begin{aligned}
&\partial_t u=\mathcal J\bigg (\kappa u+b u_{xx}+\frac{\beta}{3}q_2(q_1+2u^2)\bigg), \\
&\partial_t q_1=2u\cdot\partial_tu,\\
&\partial_t q_2=\partial_tu\cdot q_1+u\cdot\partial_tq_1=\partial_tu\cdot q_1+2u^2\cdot\partial_tu,
\end{aligned}
\right.
\end{align}
 with the consistent initial conditions
\begin{align}
u({x},0)=u_0({x}),\ q_1({x},0)=(u_0({x}))^2,\ q_2({x},0)=u_0({x},0)q_1({x},0).
\end{align}
\end{itemize}

For simplicity, in the following construction of the energy-preserving schemes, we consider the parameter $p=2$ for
the equation \eqref{Ros-equation}. Note that the extensions to the parameter $p>2$ are straightforward.
 \begin{thm}\label{RLW-equation-mass-th2.4} Under the periodic boundary condition, the reformulated system \eqref{2-Ros-IEQ-reformulation} conserves the mass \eqref{Ros-mass-energy}.
\end{thm}
\begin{prf} It follows from the first equation of \eqref{2-Ros-IEQ-reformulation}, together with the periodic boundary condition that
\begin{align*}
\frac{d}{dt}\mathcal{M}(t)&=(\partial_tu,1)\\
&=\Big(\mathcal J\big (\kappa u+b u_{xx}+\frac{\beta}{3}q+ \frac{2\beta}{3}u^2\big),1\Big)\\
&=-\Big(\partial_x\big (\kappa u+b u_{xx}+\frac{\beta}{3}q+ \frac{2\beta}{3}u^2\big),(1-\delta\partial_{xx}+\alpha\partial_{xxxx})^{-1}1\Big)\\
&=-\Big(\partial_x\big (\kappa u+b u_{xx}+\frac{\beta}{3}q+ \frac{2\beta}{3}u^2\big),1\Big)\\
&=0.
\end{align*}
This completes the proof.\qed
\end{prf}

 \begin{thm}\label{LI-E-NLS-thm-2.1} Under the periodic boundary condition, the reformulation \eqref{2-Ros-IEQ-reformulation} preserves the following invariants
  \begin{align}\label{Ros-Q-energy1}
   &\mathcal H_{1,1}=q-u^2=0,\\\label{Ros-Q-energy2}
   &\mathcal H(t)=\int_{\Omega}\Big(\frac{\kappa}{2}u^2-\frac{b}{2}u_x^2+\frac{\beta}{3} uq\Big)d{x},\ t\ge0.
  \end{align}
  \end{thm}
  \begin{prf}It follows from the second equation of \eqref{2-Ros-IEQ-reformulation} that
  \begin{align}
  \partial_t\mathcal H_{1,1}=\partial_tq-2u\cdot\partial_tu=0.
  \end{align}
  According to \eqref{2-IEQ-RLW-energy} and \eqref{2-Ros-IEQ-reformulation} together with the anti-adjoint property of $\mathcal J$, we have
  \begin{align*}
  \frac{d}{dt} \mathcal{H}(t)&=\int_{\Omega}\Big(\kappa u\partial_t u-b\partial_xu\partial_{xt}u+\frac{\beta}{3} (q\partial_tu+u\partial_tq) \Big)d{x}\\
  &=\int_{\Omega}\Big(\kappa u\partial_t u+b\partial_{xx}u\partial_t u+\frac{\beta}{3}(q+2u^2)\partial_tu\Big)d{x}\\
  &=\int_{\Omega}\Big(\kappa u+b\partial_{xx}u+\frac{\beta}{3}(q+2u^2)\Big)\mathcal J\bigg (\kappa u+b u_{xx}+\frac{\beta}{3}(q+2u^2)\bigg)d{x}\\
  &=0.
  \end{align*}
  This completes the proof. \qed
  \end{prf}

Next, we apply an RK method to the system \eqref{2-Ros-IEQ-reformulation} in time, together with the Fourier pseudo-spectral
method in space to give a class of fully discrete schemes for \eqref{2-Ros-IEQ-reformulation}, as follows:
\begin{shn}\label{Ros-scheme-2.3} Let $b_i,a_{ij}\ (i,j=1,\cdots,s)$ be real numbers and let $c_i=\sum_{j=1}^sa_{ij}$. For the given $(U^n,Q^n)$, an s-stage Runge-Kutta method is given by
\begin{align}\label{H-Ros-IEQ-equation1}
\left\{
\begin{aligned}
&K_i^n=\mathcal J_h\bigg (\kappa U^{ni}+b D_2U^{ni}+ \frac{\beta}{3}(Q^{ni}+2(U^{ni})^2)\bigg),\\
& U^{ni}=U^n+\tau\sum_{j=1}^{s}a_{ij}K_j^n,\  \mathcal J_h = -(I-\delta D_2+\alpha D_4)^{-1}D_1,\\
&Q^{ni}=Q^n+\tau\sum_{j=1}^{s}a_{ij}L_j^n,\ L_i^n=2U^{ni}\cdot K_i^n,\ i=1,2,\cdots,s,
\end{aligned}
\right.
\end{align}
and $(U^{n+1},Q^{n+1})$ is then updated by
\begin{align}\label{H-Ros-IEQ-equation2}
U^{n+1}=U^n+\tau\sum_{i=1}^{s}b_iK_i^n,\ Q^{n+1}=Q^n+\tau\sum_{i=1}^{s}b_{i}L_i^n.
\end{align}
\end{shn}

 \begin{thm}\label{Ros-equation-th2.7} If the coefficients of \eqref{H-Ros-IEQ-equation1} and \eqref{H-Ros-IEQ-equation2} satisfy \eqref{Ros-RK-coeff}, {\bf Scheme 2.2} preserves the following discrete invariants
  \begin{align}\label{Ros-semi-energy1}
   &H_{1,1}^{n+1}=H_{1,1}^n,\ H_{1,1}^n=Q^n-(U^n)^2,\\\label{Ros-semi-energy2}
   &\mathcal E_h^{n+1}=\mathcal E_h^n,\ \mathcal E_h^n=\frac{\kappa}{2}\|U^n\|^2+\frac{b}{2}\langle D_2U^n,U^n\rangle_h+\frac{\beta}{3} \langle U^n,Q^n\rangle_h,\ n=0,1,\cdots,.
  \end{align}
  \end{thm}
\begin{prf}It follows from \eqref{H-Ros-IEQ-equation2} that
\begin{align*}
H_{1,1}^{n+1}-H_{1,1}^{n}&=Q^{n+1}-Q^n-(U^{n+1})^2+(U^n)^2\\
&=\tau\sum_{j=1}^{s}b_{i}L_i^n-\tau\sum_{i=1}^{s}b_iK_i^n\cdot U^n-\tau\sum_{i=1}^{s}b_iU^n\cdot K_i^n-\tau^2\sum_{i,j=1}^{s}b_ib_jK_i^n\cdot K_j^n.
\end{align*}
Using the equality $U^n=U^{ni}-\tau\sum_{j=1}^sa_{ij}K_j^n$ and $L_i^n=2U^{ni}\cdot K_i^n$, together with \eqref{Ros-RK-coeff}, we can obtain from the above equation
\begin{align}\label{Ros-equation-2.34}
H_{1,1}^{n+1}-H_{1,1}^{n}&=\tau\sum_{j=1}^{s}b_{i}L_i^n-2\tau\sum_{i=1}^{s}b_iK_i^n\cdot U^{ni}-\tau^2\sum_{i,j=1}^{s}(b_ib_j-b_ia_{ij}-b_ja_{ji})K_i^n\cdot K_j^n\nonumber\\
&={\bf 0}.
\end{align}
With noting the initial condition $Q^0-(U^0)^2=0$, we obtain \eqref{Ros-semi-energy1} from \eqref{Ros-equation-2.34}.

According to \eqref{H-Ros-IEQ-equation1} and \eqref{H-Ros-IEQ-equation2}, we have
\begin{align*}
\mathcal E_h^{n+1}- \mathcal E_h^n&=\frac{\kappa}{2}\langle U^{n+1},U^{n+1}\rangle_h+\frac{b}{2}\langle D_2U^{n+1},U^{n+1}\rangle_h+\frac{\beta}{3} \langle U^{n+1},Q^{n+1}\rangle_h\nonumber\\
&~~~-\frac{\kappa}{2}\langle U^{n},U^{n}\rangle_h-\frac{b}{2}\langle D_2U^{n},U^{n}\rangle_h-\frac{\beta}{3} \langle U^{n},Q^{n}\rangle_h\nonumber\\
&=\frac{\tau\kappa}{2}\sum_{i=1}^{s}b_i\Big(\langle U^n,K_i^n\rangle_h+\langle K_i^n,U^n\rangle_h\Big)+\frac{\tau^2\kappa}{2}\sum_{i,j=1}^{s}b_ib_j\langle K_i^n,K_j^n\rangle_h\nonumber\\
&~~~+\frac{\tau b}{2}\sum_{i=1}^{s}b_i\Big(\langle D_2U^n,K_i^n\rangle_h+ \langle D_2K_i^n,U^n\rangle_h\Big)+\frac{\tau^2b}{2}\sum_{i,j=1}^{s}b_ib_j\langle D_2K_i^n,K_j^n\rangle_h\nonumber\\
&~~~+\frac{\tau \beta}{3}\sum_{i=1}^{s}b_i\Big(\langle U^n,L_i^n\rangle_h+ \langle Q^n,K_i^n\rangle_h\Big)+\frac{\tau^2\beta}{3}\sum_{i,j=1}^{s}b_ib_j\langle K_i^n,L_j^n\rangle_h\nonumber\\
&=\tau\sum_{i=1}^{s}b_i\Big(\kappa\langle U^{ni},K_i^n\rangle_h+b\langle D_2U^{ni},K_i^n\rangle_h+\frac{\beta}{3}\langle U^{ni},L_i^n\rangle_h+\frac{\beta}{3}\langle Q^{ni},K_i^n\rangle_h\Big)\nonumber\\
&~~~+\frac{\tau^2}{2}\kappa\sum_{i,j=1}^{s}(b_ib_j-b_ia_{ij}-b_ja_{ji})\langle K_i^n,K_j^n\rangle_h\nonumber\\
&~~~+\frac{\tau^2}{2}\beta\sum_{i,j=1}^{s}(b_ib_j-b_ia_{ij}-b_ja_{ji})\langle D_2K_i^n,K_j^n\rangle_h\nonumber\\
&~~~+\frac{\tau^2}{3}\beta\sum_{i,j=1}^{s}(b_ib_j-b_ia_{ij}-b_ja_{ji})\langle K_i^n,L_j^n\rangle_h.
\end{align*}
   With the condition \eqref{Ros-RK-coeff} and the antisymmetry of $\mathcal J_h$, we have
    \begin{align*}
    \kappa\langle &U^{ni},K_i^n\rangle_h+b\langle D_2U^{ni},K_i^n\rangle_h+\frac{\beta}{3}\langle U^{ni},L_i^n\rangle_h+\frac{\beta}{3}\langle Q^{ni},K_i^n\rangle_h\\
    &=\kappa\langle U^{ni},K_i^n\rangle_h+b\langle D_2U^{ni},K_i^n\rangle_h+\frac{2\beta}{3}\langle (U^{ni})^2,K_i^n\rangle_h+\frac{\beta}{3}\langle Q^{ni},K_i^n\rangle_h\\
    &=\langle\kappa U^{ni}+bD_2U^{ni}+\frac{2\beta}{3}(U^{ni})^2+\frac{\beta}{3}Q^{ni},K_i^n\rangle_h\\
    &=\Big\langle K_i^n,\mathcal J_hK_i^n\Big\rangle_h\\
    &=0,
    \end{align*}
which implies $ \mathcal E_h^{n+1}=\mathcal E_h^n$.
    This completes the proof. \qed
   \end{prf}
  \begin{rmk} It follows from \eqref{Ros-semi-energy1} and \eqref{Ros-semi-energy2} that under the condition \eqref{Ros-RK-coeff}, the proposed {\bf Scheme \ref{Ros-scheme-2.3}} can conserve the following discrete Hamiltonian energy
  \begin{align}\label{Ros-H-energy}
   H_h^n=\frac{\kappa}{2}\|U^n\|_h^2+\frac{b}{2}\langle D_2U^n,U^n\rangle_h+\frac{\beta}{3} \langle (U^n)^3,{\bf 1} \rangle_h,\ n=0,1,2,\cdots,.
  \end{align}

  \end{rmk}

{\begin{thm}\label{Ros-Theorem 2.8-mass} For any RK method, {\bf Scheme 2.2} conserves the following discrete mass
\begin{align}\label{semi-discrete-mass}
M_h^{n+1}=M_h^n,\ {M}_h^n=\langle U^n,{\bf 1}\rangle_h,\ n=0,1,2,\cdots,.
\end{align}
\end{thm}
\begin{prf} It follows from the first equation of \eqref{H-Ros-IEQ-equation2} that
\begin{align}
{M}_h^{n+1}={M}_h^n+\tau\sum_{i=1}^{s}b_i\langle K_i^n,{\bf 1}\rangle_h.
\end{align}
Then, we can deduce from the first equation of \eqref{H-Ros-IEQ-equation1} that
  \begin{align*}
  \langle K_i^n,{\bf 1}\rangle_h&=-\Big\langle D_1\big (\kappa U^{ni}+b D_2U^{ni}+ \frac{\beta}{3}(Q^{ni}+2(U^{ni})^2)\big),\mathcal A_h^{-1}{\bf 1}\Big\rangle_h\\
  &=-\Big\langle D_1\big (\kappa U^{ni}+b D_2U^{ni}+ \frac{\beta}{3}(Q^{ni}+2(U^{ni})^2)\big),{\bf 1}\Big\rangle_h\\
  &=0.
  \end{align*}
  This completes the proof.\qed
\end{prf}}
\begin{rmk} Here, we introduce two notations, as follows:
\begin{itemize} \item 4th-EPS: using the 4th-Gauss method to {\bf Scheme 2.2};
\item 6th-EPS: using the 6th-Gauss method to {\bf Scheme 2.2}.
\end{itemize}
\end{rmk}

\section{An efficient implementation for the proposed schemes}\label{Sec:LI-EP:3}
{As it turns out, there is no existing explicit RK methods that satisfy the condition \eqref{Ros-RK-coeff} (see Proposition 7.1.1 in \cite{FQ10}).} Thus,
in this section, motivated by \cite{CJW2021CPC,ZJ2202Arix}, we propose an efficient fixed-point iteration solver for the nonlinear equations of the proposed
schemes. For simplicity, we consider the 4th-MP scheme where the RK coefficient is given in Table \ref{Gauss-cllocation-method}. Note that the extensions to $s>2$ are straightforward.

For a given $U^n$, according to {\bf Scheme 2.1}, the 4th-MP scheme is equivalent to
\begin{align}\label{Ros-momentum-scheme-fast1}
&\mathcal A_hK_1^n=-\big[\kappa D_1U^{n1}+b D_3U^{n1}+\frac{p\beta}{p+1}((U^{n1})^{p-1}\cdot D_1U^{n1}+D_1((U^{n1})^{p}))\big],\\\label{Ros-momentum-scheme-fast2}
&\mathcal A_hK_2^n=-\big[\kappa D_1U^{n2}+b D_3U^{n2}+\frac{p\beta}{p+1}((U^{n2})^{p-1}\cdot D_1U^{n2}+D_1((U^{n2})^{p}))\big],\\\label{Ros-momentum-scheme-fast3}
& U^{n1}=U^n+\tau a_{11}K_1^n+\tau a_{12}K_2^n,\ U^{n2}=U^n+\tau a_{21}K_1^n+\tau a_{22}K_2^n,\\\label{Ros-momentum-scheme-fast4}
&U^{n+1}=U^n+\tau b_1K_1^n+\tau b_2K_2^n,
\end{align}
It follows from \eqref{Ros-momentum-scheme-fast1}-\eqref{Ros-momentum-scheme-fast3} that
\begin{align}\label{Ros-momentum-scheme-fast5}
&\mathcal A_hK_1^n=-\big[\tau\kappa a_{11}D_1K_1^n+\tau\kappa a_{12}D_1K_2^n+\tau ba_{11}D_3K_1^n+\tau b a_{12}D_3K_2^n+F_1\big],\\\label{Ros-momentum-scheme-fast6}
&\mathcal A_hK_2^n=-\big[\tau\kappa a_{21}D_1K_1^n+\tau\kappa a_{22}D_1K_2^n+\tau ba_{21}D_3K_1^n+\tau b a_{22}D_3K_2^n+F_2\big],
\end{align}
where
\begin{align*}
F_1&=\kappa D_1U^n+b D_3U^{n}+\frac{p\beta}{p+1}((U^{n1})^{p-1}\cdot D_1U^{n1}+D_1((U^{n1})^{p})),\\
F_2&=\kappa D_1U^n+b D_3U^{n}+\frac{p\beta}{p+1}((U^{n2})^{p-1}\cdot D_1U^{n2}+D_1((U^{n2})^{p})).
\end{align*}
Using Remark \ref{Ros-remark-2.1}, we can obtain from \eqref{Ros-momentum-scheme-fast5}-\eqref{Ros-momentum-scheme-fast6}
\begin{align}\label{Ros-momentum-scheme-fast7}
& A_h\mathbb{K}_1^n=-\big[\tau\kappa a_{11}\tilde{\Lambda}\mathbb{K}_1^n+\tau\kappa a_{12}\tilde{\Lambda}\mathbb K_2^n+\tau ba_{11}\tilde{\Lambda}^3\mathbb K_1^n+\tau b a_{12}\tilde{\Lambda}^3\mathbb K_2^n+\mathbb F_1\big],\\\label{Ros-momentum-scheme-fast8}
& A_h\mathbb{K}_2^n=-\big[\tau\kappa a_{21}\tilde{\Lambda}\mathbb K_1^n+\tau\kappa a_{22}\tilde{\Lambda}\mathbb K_2^n+\tau ba_{21}\tilde{\Lambda}^3\mathbb K_1^n+\tau b a_{22}\tilde{\Lambda}^3\mathbb K_2^n+\mathbb F_2\big],
\end{align}
where $ A_h=I-\delta \Lambda^2+\alpha \Lambda^4$ and $\mathbb W=\mathcal F_N W$.

For the nonlinear algebraic equations as above, we apply the following fixed-point iteration strategy, for $l=0,1,2,\cdots,M$
\begin{align*}
&\left[
  \begin{array}{cc}\vspace{2mm}
     A_h+\tau\kappa a_{11}\tilde{\Lambda}+ \tau ba_{11}\tilde{\Lambda}^3& \tau\kappa a_{12}\tilde{\Lambda}+\tau b a_{12}\tilde{\Lambda}^3\\
    \tau\kappa a_{21}\tilde{\Lambda}+\tau b a_{21}\tilde{\Lambda}^3&  A_h+\tau\kappa a_{22}\tilde{\Lambda}+ \tau ba_{22}\tilde{\Lambda}^3\\
  \end{array}
\right]\left[
  \begin{array}{c} \vspace{2mm}
  \mathbb{K}_1^{n,l+1}\\
  \mathbb{K}_2^{n,l+1}
\end{array}
  \right]=\left[
  \begin{array}{c}\vspace{2mm}
  -\mathbb F_1^{l}\\
  -\mathbb F_2^{l}
\end{array}
  \right],
\end{align*}
which implies that
\begin{align*}
\left[
  \begin{array}{cc}\vspace{2mm}
    (A_h)_j+\tau\kappa a_{11}\tilde{\Lambda}_j+ \tau ba_{11}\tilde{\Lambda}_j^3& \tau\kappa a_{12}\tilde{\Lambda}_j+\tau b a_{12}\tilde{\Lambda}_j^3\\
    \tau\kappa a_{21}\tilde{\Lambda}_j+\tau b a_{21}\tilde{\Lambda}_j^3& (A_h)_j+\tau\kappa a_{22}\tilde{\Lambda}_j+ \tau ba_{22}\tilde{\Lambda}_j^3\\
  \end{array}
\right]&\left[
  \begin{array}{c} \vspace{2mm}
  (\mathbb{K}_1^{n,l+1})_j\\
  (\mathbb{K}_2^{n,l+1})_j
\end{array}
  \right]\\
  &=\left[
  \begin{array}{c}\vspace{2mm}
  -(\mathbb F_1^{l})_j\\
  -(\mathbb F_2^{l})_j
\end{array}
  \right],
\end{align*}
where $j=0,1,\cdots,N-1$.

Solving the above equations, we obtain $\mathbb{K}_1^{n,l+1}$ and $\mathbb{K}_2^{n,l+1}$. Then ${K}_1^{n,l+1}$ and ${K}_2^{n,l+1}$ are given by ${K}_1^{n,l+1}=\mathcal F_N^H\mathbb{K}_1^{n,l+1}$ and ${K}_2^{n,l+1}=\mathcal F_N^H\mathbb{K}_2^{n,l+1}$, respectively. In our computations, we take the iterative initial value $K_1^{n,0}=U^n$ and $K_2^{n,0}=U^n$.
The iteration terminates when  the number of maximum iterative step  $M = 30$  is reached or the infinity norm of the error
between two adjacent iterative steps is less than $10^{-14}$, that is
\begin{align*}
 \mathop{\rm max}\limits_{l\leqslant i\leqslant 2} \big \{  \lVert K_i^{n,l+1}-K_i^{n,l} \rVert_{\infty,h} \big \} < 10^{-14}.
\end{align*}
Finally, we have $U^{n+1}=U^n+\tau b_1K_1^{n,M}+\tau b_2K_2^{n,M}$.

\begin{rmk}Similarly, the efficient iteration solver for the resulting nonlinear equations of the 4th-EP scheme (see {\bf Scheme 2.2}) is given by, as follows:
\begin{align*}
&\Bigg[
  \begin{array}{cc}\vspace{2mm}
    I-\tau\kappa a_{11}J_h-\tau ba_{11}J_h\Lambda^2& -\tau\kappa a_{12}J_h-\tau ba_{12}J_h\Lambda^2\\
    -\tau\kappa a_{21}J_h-\tau ba_{21}J_h\Lambda^2& I-\tau\kappa a_{22}J_h-\tau ba_{22}J_h\Lambda^2\\
  \end{array}
\Bigg]\Bigg[
  \begin{array}{c} \vspace{2mm}
  \mathbb{K}_1^{n,l+1}\\
  \mathbb{K}_2^{n,l+1}
\end{array}
 \Bigg]=\Bigg[
  \begin{array}{c}\vspace{2mm}
  \mathbb F_1^{l}\\
  \mathbb F_2^{l}
\end{array}
 \Bigg],
\end{align*}
where
\begin{align*}
F_1&= \mathcal J_h\Big(\kappa U^{n}+b D_2U^{n}+\frac{\beta}{3}\big(Q^{n1}+2(U^{n1})^2\big)\Big),\\
F_2&= \mathcal J_h\Big(\kappa U^{n}+b D_2U^{n}+\frac{\beta}{3}\big(Q^{n2}+2(U^{n2})^2\big)\Big),\\
U^{ni}&=U^n+\tau a_{i1}K_1^n+\tau a_{i2}K_2^n,\  J_h = -(I-\delta \Lambda^2+\alpha \Lambda^4)^{-1}\tilde{\Lambda},\\
Q^{ni}&=(U^n)^2+\tau\sum_{j=1}^{s}a_{ij}L_j^n,\ L_i^n=2U^{ni}\cdot K_i^n,\ i=1,2,
\end{align*}
In particular, the iteration equations as stated above  can be rewritten into the following subsystems
\begin{align*}
\left[
  \begin{array}{cc}\vspace{2mm}
    1-\tau\kappa a_{11}(J_h)_j-\tau ba_{11}(J_h)_j\Lambda_j^2& -\tau\kappa a_{12}(J_h)_j-\tau ba_{12}(J_h)_j\Lambda_j^2\\
    -\tau\kappa a_{21}(J_h)_j-\tau ba_{21}(J_h)_j\Lambda_j^2& 1-\tau\kappa a_{22}(J_h)_j-\tau ba_{22}(J_h)_j\Lambda_j^2\\
  \end{array}
\right]&\left[
  \begin{array}{c} \vspace{2mm}
  (\mathbb{K}_1^{n,l+1})_j\\
  (\mathbb{K}_2^{n,l+1})_j
\end{array}
  \right]\\
  &=\left[
  \begin{array}{c}\vspace{2mm}
  (\mathbb F_1^{l})_j\\
  (\mathbb F_2^{l})_j
\end{array}
  \right],
\end{align*}
where $j=0,1,2,\cdots,N-1$. Then, by the similar procedure as the 4th-MP scheme, we have $U^{n+1}$.
\end{rmk}

 \section{Numerical results}\label{Ros-section-3}
In this section, we devote to the numerical performances in the accuracy and
invariants-preservation of the proposed schemes for the equation \eqref{Ros-equation}. For brevity, in the rest
of this paper, 4th-MPS, 4th-EPS, 6th-MPS and 6th-EPS are
only used for demonstration purposes. In order to quantify the numerical solution, we introduce the $l^2$-and $l^{\infty}$-norm error functions, respectively,
\begin{align*}
&e_{\infty}(t_n)=\|u(\cdot,t_n)-u^n\|_{h,\infty},\ e_{2}(t_n)=\|u(\cdot,t_n)-u^n\|_{h}.
\end{align*}
Furthermore, we also investigate the residuals on the mass, momentum and Hamiltonian energy, defined respectively, as
\begin{align*}
Erorr_1^n=\big|M_h^n-M_h^0\big|,\ Erorr_2^n=\big|I_h^n-I_h^0\big|,\ Erorr_3^n=\big|H_h^n-H_h^0\big|.
\end{align*}
All simulations are performed on a Win10 machine with Intel Core i7 and 32GB
using MATLAB R2015b.


\subsection{Rosenau-RLW equation}

As the parameters $\kappa=\delta=\alpha=\beta=1$ and $b=0$ are chosen, the equation \eqref{Ros-equation} reduces to the generalized Rosenau-RLW equation, which has an exact solution \cite{CSWAMC-2015,PZAMM2012}
\begin{align}
u(x,t)=A\text{sech}^{\frac{4}{p-1}}\big(B(x-Ct-x_0)\big),
\end{align}
where
\begin{align*}
A = \exp\Big(\frac{\ln\frac{(p+3)(3p+1)(p+1)}{2(p^2+3)(p^2+4p+7)}}{p-1}\Big),\ B = \frac{p-1}{\sqrt{4p^2+8p+20}},\  C = \frac{p^4+4p^3+14p^2+20p+25}{p^4+4p^3+10p^2+12p+21},
\end{align*}
and $x_0$ represents the initial phase of the solition and the periodic boundary condition is considered.

We first verify the convergence order in time for the proposed schemes. Let us set the computational domain $\Omega=[-200,200]$ and we take the initial value
\begin{align*}
u(x,0)=A\text{sech}^{\frac{4}{p-1}}\big(B(x-x_0)\big),\ x\in\Omega,
\end{align*}
where the initial phase $x_0=0$.

\begin{figure}[H]
\centering\begin{minipage}[t]{60mm}
\includegraphics[width=60mm]{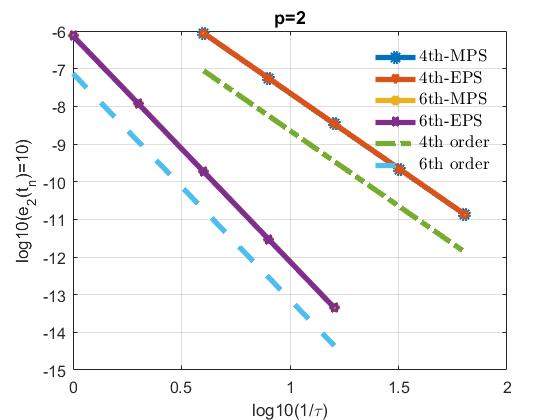}
\end{minipage}
\begin{minipage}[t]{60mm}
\includegraphics[width=60mm]{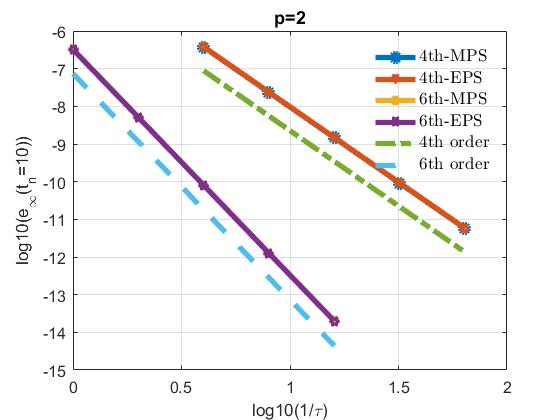}
\end{minipage}
\begin{minipage}[t]{60mm}
\includegraphics[width=60mm]{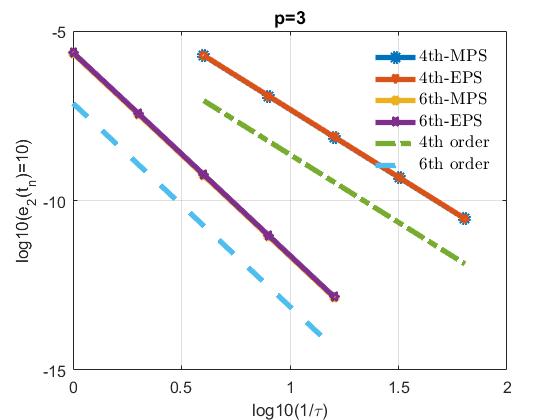}
\end{minipage}
\centering\begin{minipage}[t]{60mm}
\includegraphics[width=60mm]{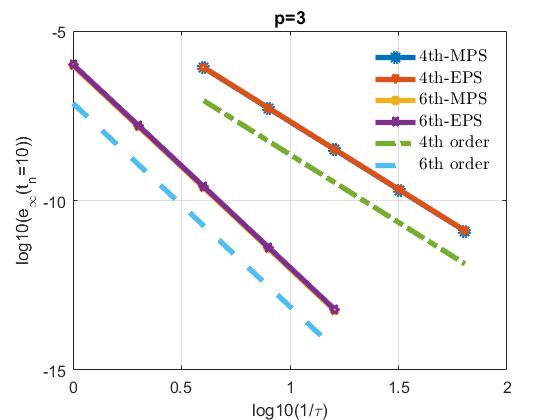}
\end{minipage}
\begin{minipage}[t]{60mm}
\includegraphics[width=60mm]{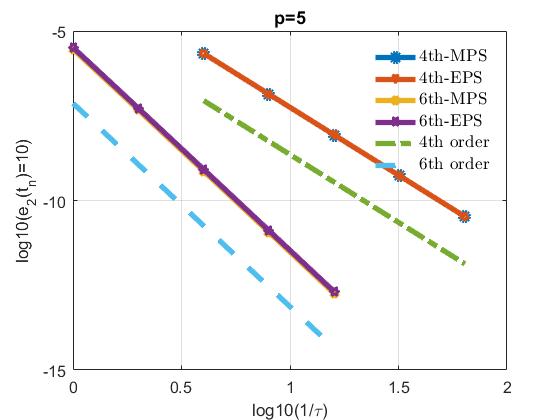}
\end{minipage}
\begin{minipage}[t]{60mm}
\includegraphics[width=60mm]{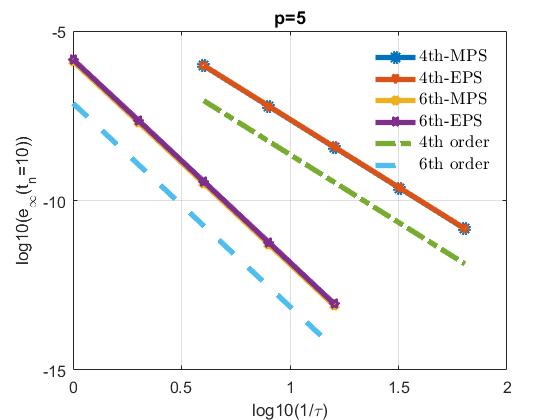}
\end{minipage}
 \caption{{The $l^2$ and $l^{\infty}$-norm errors vs. the time step sizes provided by the proposed 4th-MPS, 4th-EPS, 6th-MPS and 6th-EPS for the Rosenau-RLW equation with different  parameters $p$, respectively.}}\label{R-RLW1d-error}
\end{figure}

The temporal convergence tests are investigated by fixing the Fourier node $2048$. We compute the
numerical solutions at $t=10$ using 4th-MPS and 4th-EPS with various time step sizes $\tau=2^{-k},\ k=2,3,4,5,6$ as well as 6th-MPS and 6th-EPS with various time step sizes $\tau=2^{-k},\ k=0,1,2,3,4$, respectively. The relation between the $l^2$ and $l^{\infty}$-norm errors and the time step size is summarized in Figure \ref{R-RLW1d-error},
where the up picture corresponds to the parameter $p=2$, the middle one corresponds to the parameter $p=3$ and the down one is for the parameter $p=5$. The fourth-order temporal accuracy for 4th-MPS and 4th-EPS and the sixth-order for 6th-MPS and 6th-EPS are observed, respectively, for different parameters $p$ as desired.

\begin{figure}[H]
\centering\begin{minipage}[t]{45mm}
\includegraphics[width=45mm]{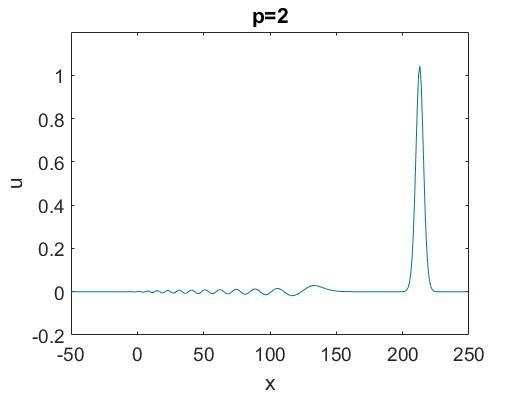}
\end{minipage}
\begin{minipage}[t]{45mm}
\includegraphics[width=45mm]{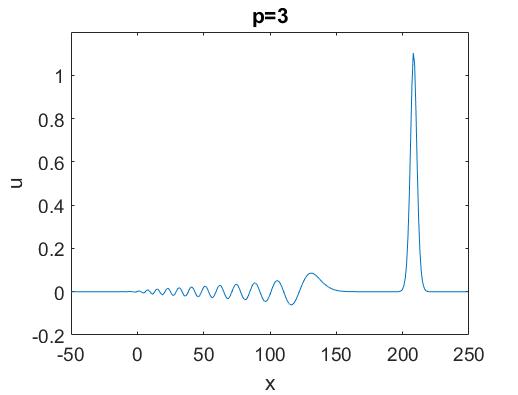}
\end{minipage}
\begin{minipage}[t]{45mm}
\includegraphics[width=45mm]{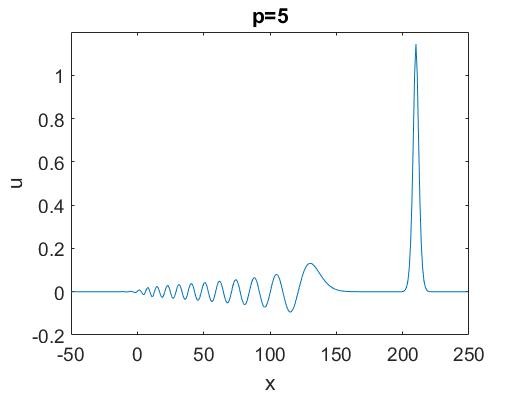}
\end{minipage}
 \caption{{The profile of $u$ provided by 4th-MPS with different parameters $p$ at $t=100$, where the uniform spatial and time step size is chosen as $\tau=h = 0.1$ for the Rosenau-RLW equation.}}\label{R-RLW1d-soliton}
\end{figure}

\begin{figure}[H]
\centering\begin{minipage}[t]{45mm}
\includegraphics[width=45mm]{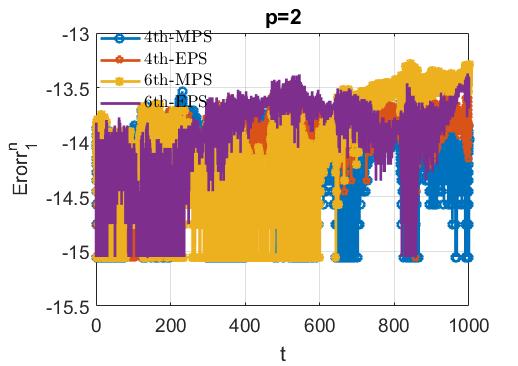}
\end{minipage}
\begin{minipage}[t]{45mm}
\includegraphics[width=45mm]{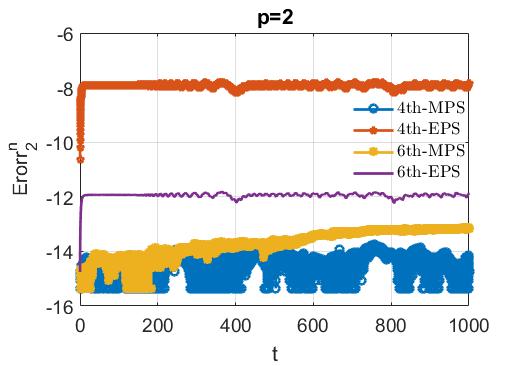}
\end{minipage}
\begin{minipage}[t]{45mm}
\includegraphics[width=45mm]{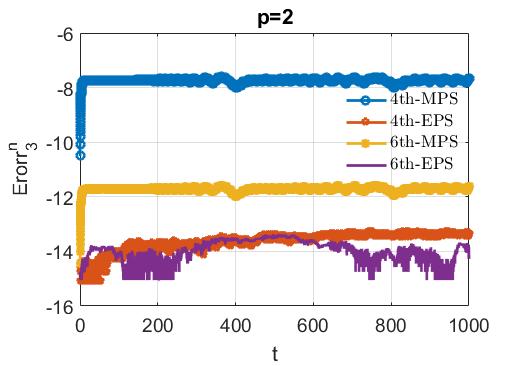}
\end{minipage}
\centering\begin{minipage}[t]{45mm}
\includegraphics[width=45mm]{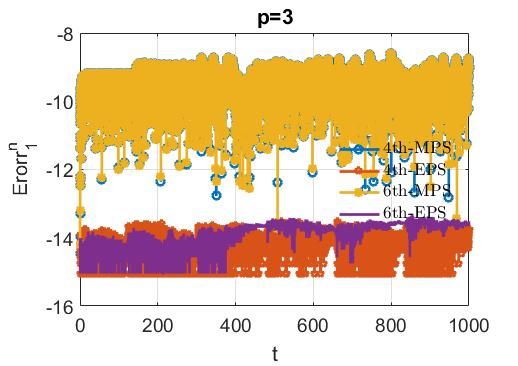}
\end{minipage}
\begin{minipage}[t]{45mm}
\includegraphics[width=45mm]{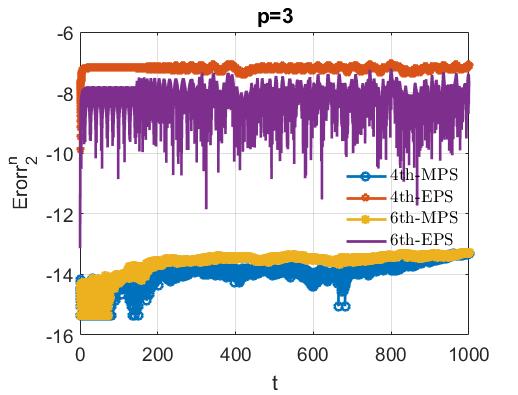}
\end{minipage}
\begin{minipage}[t]{45mm}
\includegraphics[width=45mm]{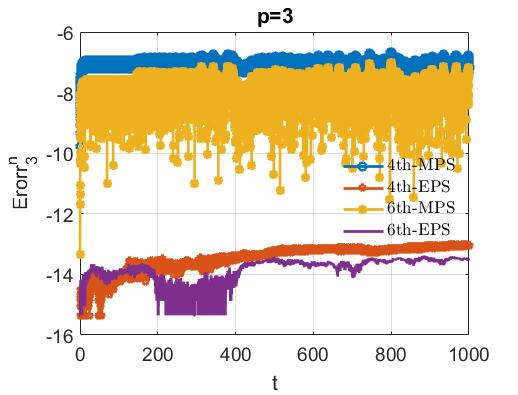}
\end{minipage}
\centering\begin{minipage}[t]{45mm}
\includegraphics[width=45mm]{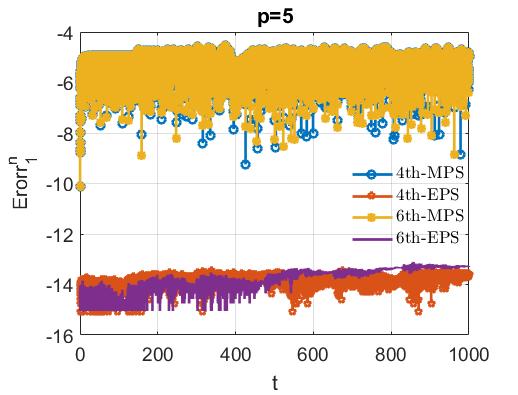}
\end{minipage}
\begin{minipage}[t]{45mm}
\includegraphics[width=45mm]{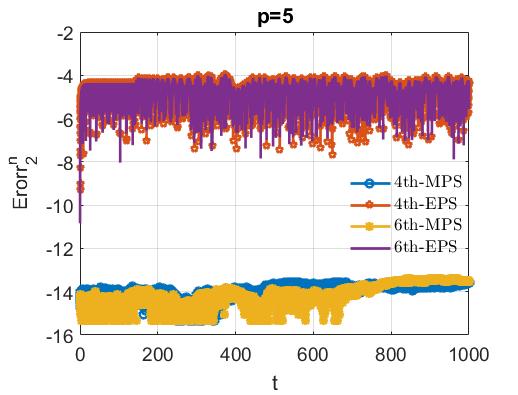}
\end{minipage}
\begin{minipage}[t]{45mm}
\includegraphics[width=45mm]{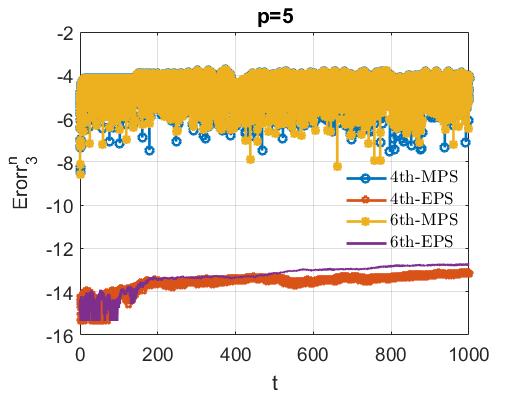}
\end{minipage}
 \caption{{The residuals on the discrete mass, momentum and Hamiltonian energy from $t=0$ to $t=1000$ provided by the proposed 4th-MPS, 4th-EPS, 6th-MPS and 6th-EPS for the Rosenau-RLW equation with different parameters $p$, respectively.}}\label{R-RLW-residuals}
\end{figure}

Then, we choose the initial value, as follows:
\begin{align}
u(x,0) = \exp(-0.05(x-40)^2),\ x\in \Omega,
\end{align}
where the computational domain $\Omega=[-50,250]$ with the periodic boundary condition. We take the uniform spatial size $h=1$ and time step size $\tau=0.1$, respectively, and compute the profile of $u$ at $T=100$ by choosing different parameters $p$ as well as the residuals on the discrete mass, momentum and Hamiltonian energy by using 4th-MPS, 4th-EPS, 6th-MPS and 6th-EPS, respectively.
Figure \ref{R-RLW1d-soliton} shows the profile of $u$ provided by 4th-MPS at $t=100$, where the left picture corresponds to the parameter $p=2$, the middle one corresponds to the parameter $p=3$ and the right one is for the parameter $p=5$. The influence of the parameter $p$ on the dispersion wave propagation is observed as shown in \cite{DLN2021-MCS}. We should note that the
profiles of $u$ at $t=100$ computed by using other schemes are similar to Figure \ref{R-RLW1d-soliton}, thus for brevity, we omit them. Figure \ref{R-RLW-residuals} shows the evolutions of the residuals on the discrete mass, momentum and Hamiltonian energy of numerical solutions computed by 4th-MPS, 4th-EPS, 6th-MPS and 6th-EPS, respectively.
We can draw the following observations: (1) for the parameter $p=2$, all of proposed schemes preserve the discrete mass exactly, while the proposed momentum-preserving schemes cannot conserve the discrete mass for the parameters $p=3$ and $p=5$; (2) 4th-MPS and 6th-MPS can preserve discrete momentum up to the machine precision; (3) 4th-EPS and 6th-EPS conserve the discrete both mass and Hamiltonian energy for all parameters $p$.

\subsection{Rosenau-KdV equation}\label{R-KdV-equation-JQSZ}
As the parameters $\kappa=b=\alpha=1$ and $\delta=0$ are chosen, the equation \eqref{Ros-equation} reduces to the generalized Rosenau-KdV equation\cite{ZZ-JAM-2013}. For simplicity, we also consider the parameters $p=2$, $p=3$ and $p=5$, respectively, where the exact solution can be given by
\begin{itemize}
\item {\bf Case I}: As $\beta=\frac{1}{2}$ and $p=2$, the Rosenau-KdV equation has an exact solution \cite{HXH-AMP-2013}
\begin{align}
u(x,t)=k_{11}{\rm sech}^4(k_{12}(x-k_{13}t)).
\end{align}
where $k_{11} = -\frac{35}{24}+\frac{35}{312}\sqrt{313},\ k_{12} = \frac{1}{24}\sqrt{-26+2\sqrt{313}},\ k_{13} = \frac{1}{2}+\frac{\sqrt{313}}{26}$.
 \item {\bf Case II}: If $\beta=1$ and $p=3$, the Rosenau-KdV equation admits an exact solution  \cite{RMB-RJP-2014}
\begin{align}
u(x,t)=k_{21}{\rm sech}^2(k_{22}(x-k_{23}t)),
\end{align}
where $k_{21} = \frac{1}{4}\sqrt{-15+3\sqrt{41}},\ k_{22} = \frac{1}{4}\sqrt{\frac{-5+\sqrt{41}}{2}},\ k_{23} = \frac{1}{10}(5+\sqrt{41})$.
\item {\bf Case III}: When $\beta=1$ and $p=5$, the exact solution of the Rosenau-KdV equation is given by \cite{ZZ-JAM-2013}
\begin{align}
u(x,t)=k_{31}{\rm sech}^2(k_{32}(x-k_{33}t)),
\end{align}
where $k_{31} = \sqrt[4]{\frac{4}{15}(-5+\sqrt{34})},\ k_{32} = \frac{1}{3}\sqrt{-5+\sqrt{34}},\ k_{33} = \frac{1}{10}(5+\sqrt{34})$.
\end{itemize}

\begin{figure}[H]
\centering\begin{minipage}[t]{65mm}
\includegraphics[width=65mm]{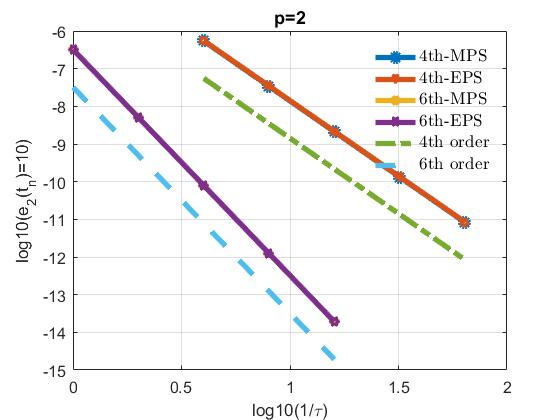}
\end{minipage}
\begin{minipage}[t]{65mm}
\includegraphics[width=65mm]{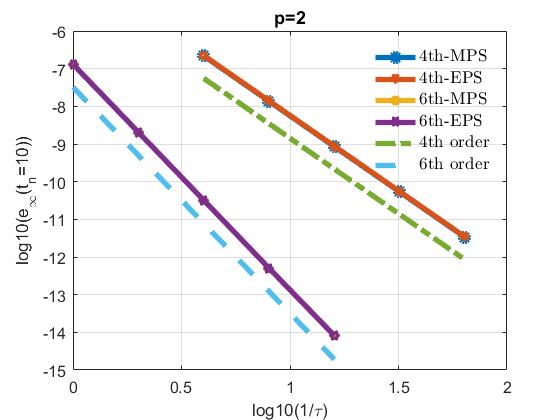}
\end{minipage}
\begin{minipage}[t]{65mm}
\includegraphics[width=65mm]{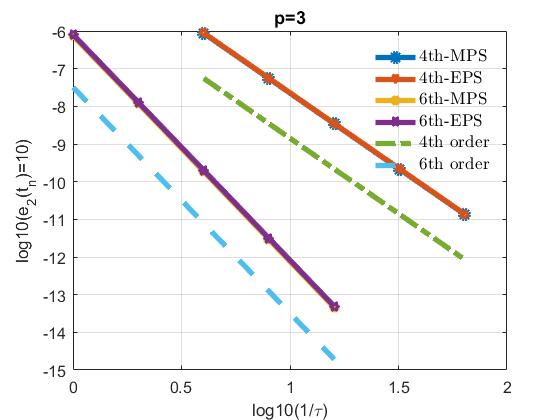}
\end{minipage}
\centering\begin{minipage}[t]{65mm}
\includegraphics[width=65mm]{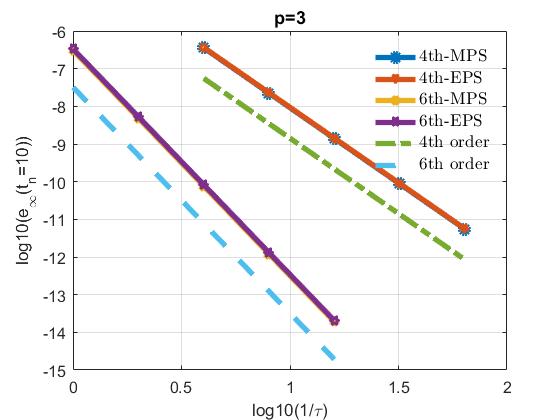}
\end{minipage}
\begin{minipage}[t]{65mm}
\includegraphics[width=65mm]{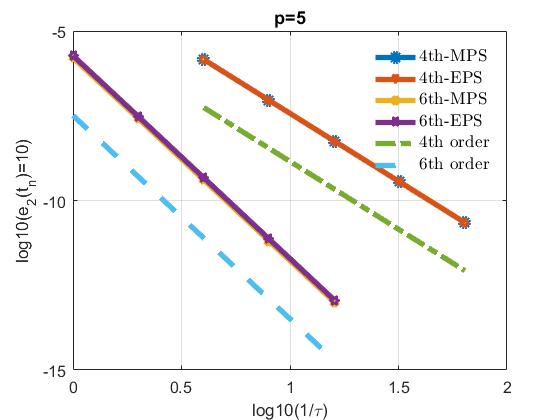}
\end{minipage}
\begin{minipage}[t]{65mm}
\includegraphics[width=65mm]{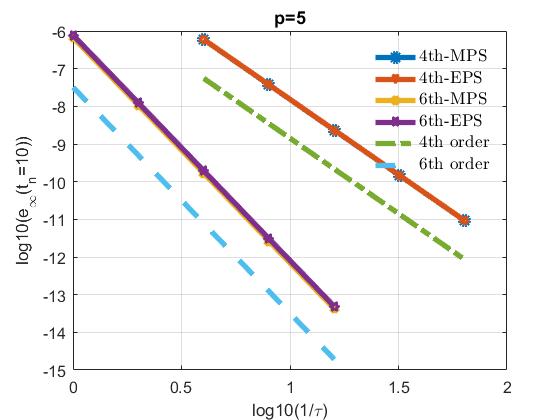}
\end{minipage}
 \caption{{The $l^2$ and $l^{\infty}$-norm errors vs. the time step size provided by the proposed 4th-MPS, 4th-EPS, 6th-MPS and 6th-EPS for the Rosenau-KdV equation with different parameters $p$, respectively.}}\label{R-RdV-error}
\end{figure}

 First of all, we verify the convergence order in time for the selected four structure-preserving schemes. Let us set the computational domain $\Omega=[-200,200]$ and we take the initial value as the exact solution at $t=0$ for the parameters $p=2$, $p=3$ and $p=5$, respectively. The temporal convergence tests are conducted by fixing the Fourier node $2048$. Figure \ref{R-RdV-error} shows the $l^2$ and $l^{\infty}$-norm errors with various time step sizes at $t=10$, where we take time step sizes $\tau=2^{-k},\ k=2,3,4,5,6$ for 4th-MPS and 4th-EPS, while for 6th-MPS and 6th-EPS, we choose time step sizes $\tau=2^{-k},\ k=0,1,2,3,4$. It is clear to see that 4th-MPS and 4th-EPS are fourth-order temporal accuracy, and 6th-MPS and 6th-EPS can achieve sixth-order accuracy in time.

\begin{figure}[H]
\centering\begin{minipage}[t]{45mm}
\includegraphics[width=45mm]{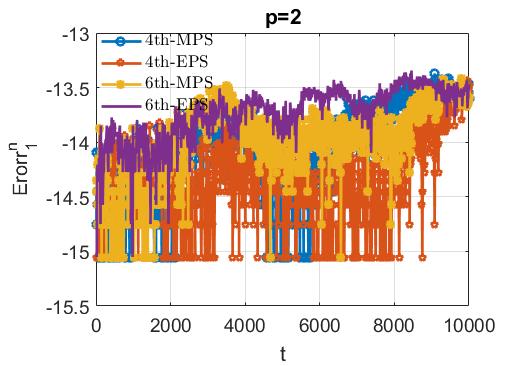}
\end{minipage}
\begin{minipage}[t]{45mm}
\includegraphics[width=45mm]{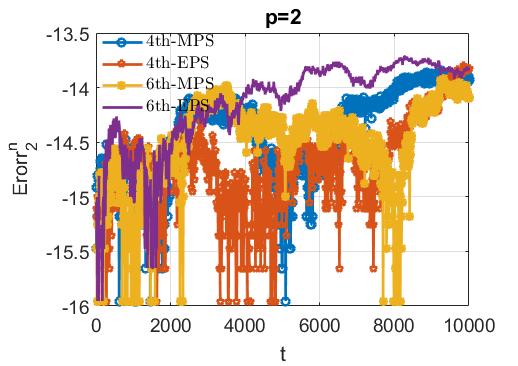}
\end{minipage}
\begin{minipage}[t]{45mm}
\includegraphics[width=45mm]{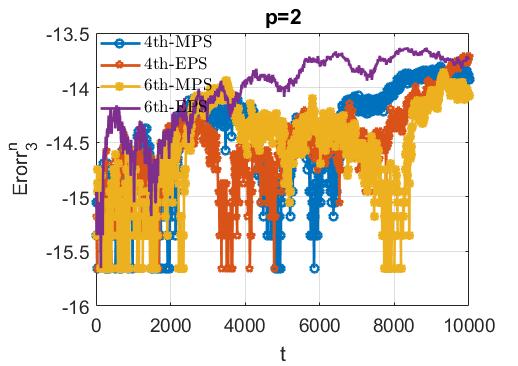}
\end{minipage}
\centering\begin{minipage}[t]{45mm}
\includegraphics[width=45mm]{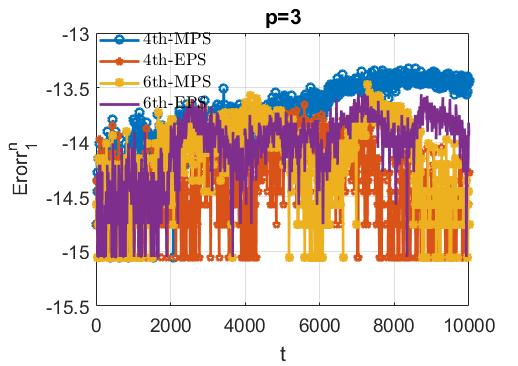}
\end{minipage}
\begin{minipage}[t]{45mm}
\includegraphics[width=45mm]{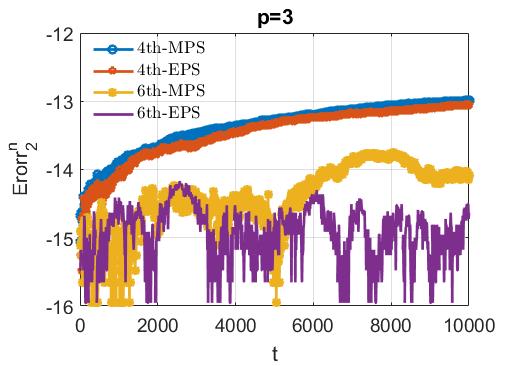}
\end{minipage}
\begin{minipage}[t]{45mm}
\includegraphics[width=45mm]{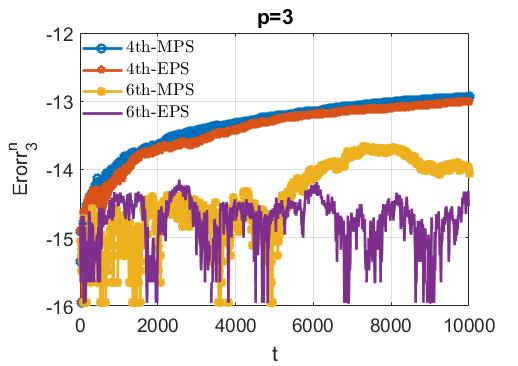}
\end{minipage}
\centering\begin{minipage}[t]{45mm}
\includegraphics[width=45mm]{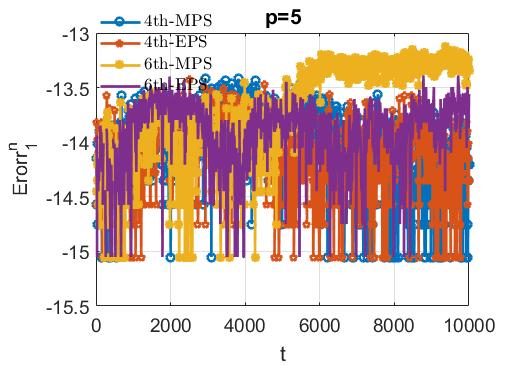}
\end{minipage}
\begin{minipage}[t]{45mm}
\includegraphics[width=45mm]{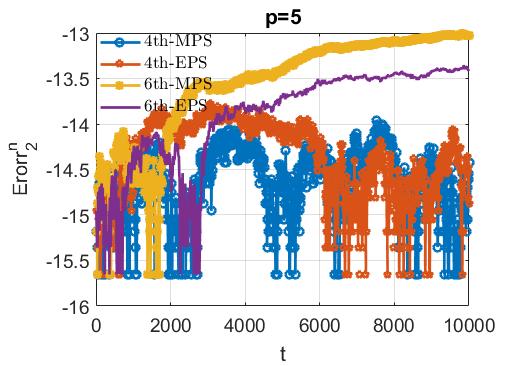}
\end{minipage}
\begin{minipage}[t]{45mm}
\includegraphics[width=45mm]{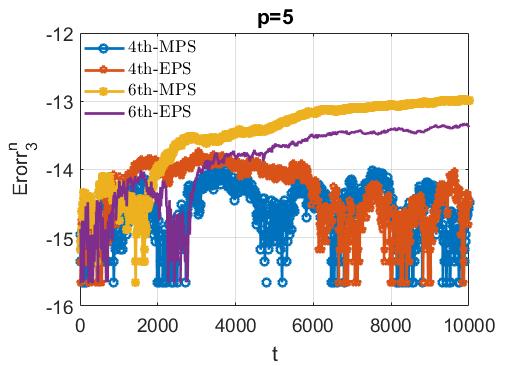}
\end{minipage}
 \caption{{The residuals on the mass, momentum and Hamiltonian energy from $t=0$ to $t=10000$ provided by the proposed 4th-MPS, 4th-EPS, 6th-MPS and 6th-EPS for the Rosenau-RdV equation with different parameters $p$, respectively.}}\label{R-KdV-residuals}
\end{figure}

Then, we set the computational domain $\Omega=[-100,100]$ and take the uniform spatial step size $h=\frac{200}{512}$ and time step size $\tau=0.1$, respectively. Figure \ref{R-KdV-residuals} shows the residuals on the discrete mass, momentum and Hamiltonian energy computed by using 4th-MPS, 4th-EPS, 6th-MPS and 6th-EPS, respectively, which is similar to Figure \ref{R-RLW-residuals}. We should note that the residuals on the discrete mass provided by 4th-MPS and 6th-MPS is also up to the machine precision because of the fine spatial mesh.

\subsection{Some comparisons}
In the following numerical experiments, we compare the proposed structure-preserving schemes with the linearized Crank-Nicolson momentum-preserving scheme (LCN-MPS) \cite{JCCWarXiv2019} and fourth-order energy-preserving schemes (i.e., YC-EPS and SC-EPS) \cite{CLZCNSNS2018} by focusing on the numerical errors in time and the computational efficiency, respectively.

Let us still consider the Rosenau-KdW equation as above. For the sake of simplicity, we only consider the parameter $p=2$ (see Case I in Section \ref{R-KdV-equation-JQSZ}).
We set the computational domain $\Omega=[-100,100]$ and the uniform spatial mesh $h=\frac{200}{1024}$. Table \ref{Tab:Rose-equation:1} shows the numerical errors and convergence order in time for the different schemes with various time steps at $t=1$. It is clear to observe that (i) LCN-MPS is second-order accurate in time and its numerical errors are largest; (ii) 4th-MPS, 4th-EPS, YC-EPS and SC-EPS are all fourth-order accurate in time, and the numerical error provided by SC-EPS is smallest, while the ones provided by YC-EPS are much larger than other fourth-order schemes.

\begin{figure}[H]
\centering\begin{minipage}[t]{65mm}
\includegraphics[width=65mm]{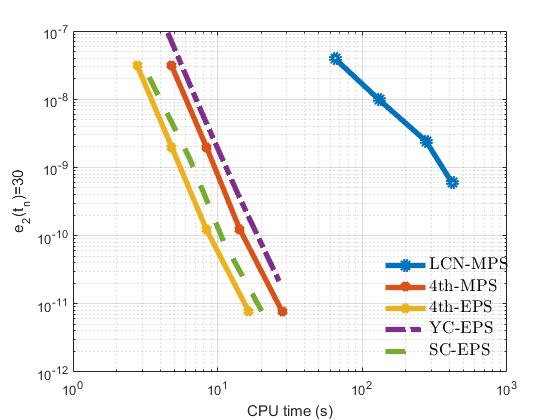}
\end{minipage}
\begin{minipage}[t]{65mm}
\includegraphics[width=65mm]{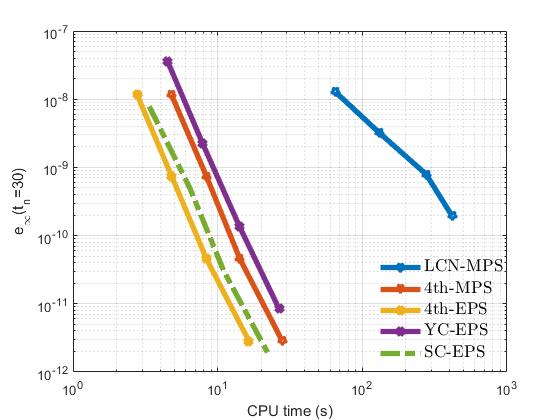}
\end{minipage}
 \caption{The $l^2$ and $l^{\infty}$-norm errors vs. the CPU time provided by LCN-MPS, 4th-MPS, 4th-EPS, YC-EPS and SC-EPS for the Rosenau-KdV equation.}\label{R-KdV-CPU}
\end{figure}

Finally, we set the computational domain $\Omega=[-300,300]$ and the Fourier node $2^{13}$, and we then investigate the global $l^2$-and $l^{\infty}$- errors of $u$ versus the CPU time using the five selected structure-preserving schemes with various time steps at $t=30$. The results are summarized in Figure \ref{R-KdV-CPU}. For a given global error, we observe that (i) the cost of the LCN-MPS is the most expensive because of the low-order accuracy in time; (ii) the cost of 4th-EPS is the cheapest; (iii) the cost of SC-EPS is much cheaper than the one provided by the 4th-MPS, and the cost of 4th-MPS is much cheaper than the one provided by the YC-EPS. We should note that as the parameter $p$ is enlarged, more QAV variables shall be introduced, thus the computational cost of the high-order energy-preserving schemes will increase.


\begin{table}[H]
\tabcolsep=9pt
\footnotesize
\renewcommand\arraystretch{1.1}
\centering
\caption{{Numerical errors and convergence order for the different schemes with various time steps at $t=1$.}}\label{Tab:Rose-equation:1}
\begin{tabular*}{\textwidth}[h]{@{\extracolsep{\fill}}c c c c c c c}\hline
{Scheme\ \ } &{$\tau$} &{$e_{2}(t_n=1)$} &{order}& {$e_{\infty}(t_n=1)$}&{order}  \\     
\hline
\multirow{4}{*}{4th-EPS}  &{$\frac{1}{10}$}& {1.585e-09}&{-} &{6.292e-10} & {-}\\[1ex]
 {}  &{$\frac{1}{20}$}& {9.905e-11}&{4.000} &{3.933e-11} & {4.000}\\[1ex] 
{}  &{$\frac{1}{40}$}& {6.191e-12}&{4.000} &{2.458e-12} &{4.000} \\[1ex]
{}  &{$\frac{1}{80}$}& {3.867e-13}&{4.001} &{1.540e-13} &{3.997} \\
 \multirow{4}{*}{4th-MPS}  &{$\frac{1}{10}$}& {1.518e-09}&{-} &{5.984e-10} & {-}\\[1ex]
 {}  &{$\frac{1}{20}$}& {9.490e-11}&{4.000} &{3.740e-11} & {4.000}\\[1ex] 
  {}  &{$\frac{1}{40}$}& {5.932e-12}&{4.000} &{ 2.338e-12} &{4.000} \\[1ex]
  {}  &{$\frac{1}{80}$}& {3.706e-13}&{4.000} &{1.462e-13} &{3.999} \\
    \multirow{4}{*}{LCN-MPS\cite{JCCWarXiv2019}}  &{$\frac{1}{10}$}& {7.299e-05}&{-} &{2.800e-05} & {-}\\[1ex]
  {}  &{$\frac{1}{20}$}& {1.813e-05}&{2.009} &{6.950e-06} & { 2.010}\\[1ex] 
   {}  &{$\frac{1}{40}$}& {4.518e-06}&{2.005} &{1.731e-06} &{2.006} \\[1ex]
    {}  &{$\frac{1}{80}$}& {1.128e-06}&{2.002} &{4.319e-07} &{2.003} \\
     \multirow{4}{*}{YC-EPS \cite{CLZCNSNS2018}}  &{$\frac{1}{10}$}& {8.024e-08}&{-} &{3.400e-08} & {-}\\[1ex]
  {}  &{$\frac{1}{20}$}& { 5.024e-09}&{3.997} &{ 2.129e-09} & {3.997}\\[1ex] 
   {}  &{$\frac{1}{40}$}& {3.142e-10}&{ 3.999} &{1.331e-10} &{3.999} \\[1ex]
    {}  &{$\frac{1}{80}$}& {1.964e-11}&{4.000} &{8.326e-12} &{4.000} \\
    \multirow{4}{*}{SC-EPS \cite{CLZCNSNS2018}}  &{$\frac{1}{10}$}& {1.081e-09}&{-} &{4.355e-10} & {-}\\[1ex]
  {}  &{$\frac{1}{20}$}& {6.761e-11}&{4.000} &{2.723e-11} & {3.999}\\[1ex] 
   {}  &{$\frac{1}{40}$}& { 4.229e-12}&{3.999} &{1.708e-12} &{3.995} \\[1ex]
    {}  &{$\frac{1}{80}$}& {2.896e-13}&{3.868} &{1.186e-13} &{3.848} \\\hline
\end{tabular*}
\end{table}

\section{Conclusions}\label{Ros-section-4}
In this paper, we propose two classes of high-order structure-preserving schemes for the generalized Rosenau-type equation \eqref{Ros-equation}. One of the schemes can conserve the discrete momentum conservation law, which is based on the use of the symplectic RK method in time and the standard Fourier pseudo-spectral method in space, respectively. Another one conservers the discrete Hamiltonian energy and mass, where the main idea is based on the combination of the QAV approach with the the symplectic RK method in time, together with the standard Fourier pseudo-spectral method in space. Extensive numerical
tests and comparisons are also addressed to verify the performance of the proposed schemes.

We conclude this paper with two main remarks. {First, we note that the construction of the energy-preserving schemes should be discussed case by case, and the proposed momentum-preserving schemes are not mass-conserving as the parameter $p>2$. Thus in practical computations, such trade-offs between two classes of schemes shall be treated carefully.} Second, as far as we know, there are some works on optimal error estimates of EQ schemes \cite{CW-2021-AMMM,LSR-2019-MP,ST-SIAM-2018,WHW-JSC-2021} and Fourier pseudo-spectral methods \cite{CFGW-2015-NMPDE,GWWC17,GTWWW-2012-Siam,ZWHWY-ANM-2017}, but  to our best knowledge, the error estimate of high-order structure-preserving Fourier pseudo-spectral schemes is still not available. Therefore, how to establish optimal error estimates for the proposed schemes will be an interesting topic for future studies. In fact, the uniformly bounded of numerical solutions in $l^{\infty}$-norm can be obtained by using the discrete momentum \eqref{Ros-equation-new-2.5} and Sobolev imbedding theorems \cite{GW12}, thus our first attempt will focus on the high-order momentum-preserving schemes, followed by the energy-preserving schemes.



 \section*{Acknowledgments}

The authors would like to express sincere gratitude to the referees for their insightful comments
and suggestions. The work is supported by the National Natural Science Foundation of China (Grant Nos. 11971481,12261097), and the Yunnan Fundamental Research Projects (Grant No. 202101AT070208), and the Natural Science Foundation of Hunan Province (Grant Nos. 2021JJ40655, 2021JJ20053).


\begin{thebibliography}{10}

\bibitem{AQJMS2022}
M.~Ahmat and J.~Qiu.
\newblock {SSP IMEX Runge-Kutta} {WENO} scheme for generalized
  {R}osenau-{KdV}-{RLW} equation.
\newblock {\em J. Math. Study.}, 55:1--21, 2022.

\bibitem{AOAA2015}
N.~Atouani and K.~Omrani.
\newblock A new conservative high-order accurate difference scheme for the
  {R}osenau equation.
\newblock {\em Appl. Anal.}, 94:2435--2455, 2015.

\bibitem{CLZCNSNS2018}
J.~Cai, H.~Liang, and C.~Zhang.
\newblock Efficient high-order structure-preserving methods for the generalized
  {R}osenau-type equation with power law nonlinearity.
\newblock {\em Commun. Nonlinear Sci. Numer. Simulat.}, 59:122--131, 2018.

\bibitem{CSWAMC-2015}
W.~Cai, Y.~Sun, and Y.~Wang.
\newblock Variational discretizations for the generalized {R}osenau-type
  equations.
\newblock {\em Appl. Math. Comput.}, 271:860--873, 2015.

\bibitem{CQ01}
J.~Chen and M.~Qin.
\newblock Multi-symplectic {F}ourier pseudospectral method for the nonlinear
  {S}chr\"{o}dinger equation.
\newblock {\em Electr. Trans. Numer. Anal.}, 12:193--204, 2001.

{\bibitem{CGHW-NMTMA-2022}
Y.~Chen, Y.~Gong, Q.~Hong, and C.~Wang.
\newblock A novel class of energy-preserving {R}unge-{K}utta methods for the
  {K}orteweg-de {V}ries equation.
\newblock {\em Numer. Math. Theor. Meth. Appl.}, 15:768--792, 2022.}

{\bibitem{CFGW-2015-NMPDE}
K.~Cheng, W.~Feng, S.~Gottlieb, and C.~Wang.
\newblock A {F}ourier pseudospectral method for the ``good" {B}oussinesq
  equation with second-order temporal accuracy.
\newblock {\em Numer. Methods Partial Differ. Equ.}, 31:202--224, 2015.}

{\bibitem{CW-2016-Siam}
K.~Cheng and C.~Wang.
\newblock Long time stability of high order multi-step numerical schemes for
  two-dimensional incompressible {N}avier-{S}tokes equations.
\newblock {\em SIAM J. Numer. Anal.}, 54:3123--3144, 2016.}

{\bibitem{CW-2021-AMMM}
Q.~Cheng and C.~Wang.
\newblock Error estimate of a second order accurate scalar auxiliary variable
  ({SAV}) scheme for the thin film epitaxial equation.
\newblock {\em Adv. Appl. Math. Mech.}, 13:1318--1354, 2021.}

\bibitem{ChungAA98}
S.~K. Chung.
\newblock Finite difference approximate solutions for the {R}osenau equation.
\newblock {\em Appl. Anal.}, 69:149--156, 1998.

\bibitem{CJW2021CPC}
J.~Cui, Y.~Wang, and C.~Jiang.
\newblock Arbitrarily high-order structure-preserving schemes for the
  {Gross-Pitaevskii} equation with angular momentum rotation.
\newblock {\em Comput. Phys. Commun.}, 261:107767, 2021.

\bibitem{DLN2021-MCS}
Y.~I. Dimitrienko, S.~Li, and Y.~Niu.
\newblock Study on the dynamics of a nonlinear dispersion model in both {1D and
  2D} based on the fourth-order compact conservative difference scheme.
\newblock {\em Math. Comput. Simulat.}, 182:661--689, 2021.

\bibitem{FQ10}
K.~Feng and M.~Qin.
\newblock {\em Symplectic Geometric Algorithms for {H}amiltonian Systems}.
\newblock Springer and Zhejiang Science and Technology Publishing House,
  Heidelberg, Hangzhou, 2010.

\bibitem{FM2011}
D.~Furihata and T.~Matsuo.
\newblock {\em Discrete Variational Derivative Method: A Structure-Preserving
  Numerical Method for Partial Differential Equations}.
\newblock Chapman \& Hall/CRC, Boca Raton, 2011.

\bibitem{GCW14}
Y.~Gong, J.~Cai, and Y.~Wang.
\newblock Multi-symplectic {F}ourier pseudospectral method for the {K}awahara
  equation.
\newblock {\em Commun. Comput. Phys.}, 16:35--55, 2014.

\bibitem{GHWW2022}
Y.~Gong, Q.~Hong, C.~Wang, and Y.~Wang.
\newblock Quadratic auxiliary variable {Runge-Kutta} methods for the
  {C}amassa-{H}olm equation.
\newblock {\em Adv. Appl. Math. Mech.}, 2023, 10.4208/aamm.OA-2022-0188.

\bibitem{GWWC17}
Y.~Gong, Q.~Wang, Y.~Wang, and J.~Cai.
\newblock A conservative {F}ourier pseudo-spectral method for the nonlinear
  {S}chr\"{o}dinger equation.
\newblock {\em J. Comput. Phys.}, 328:354--370, 2017.

\bibitem{GZW2020jcp}
Y.~Gong, J.~Zhao, and Q.~Wang.
\newblock Arbitrarily high-order linear energy stable schemes for gradient flow
  models.
\newblock {\em J. Comput. Phys.}, 419:109610, 2020.

{\bibitem{GTWWW-2012-Siam}
S.~Gottlieb, F.~Tone, C.~Wang, X.~Wang, and D.~Wirosoetisno.
\newblock Long time stability of a classical efficient scheme for
  two-dimensional {N}avier-{S}tokes equations.
\newblock {\em SIAM J. Numer. Anal.}, 50:126--150, 2012.}

{\bibitem{GW12}
S.~Gottlieb and C.~Wang.
\newblock Stability and convergence analysis of fully discrete {F}ourier
  collocation spectral method for 3-{D} viscous {B}urgers' equation.
\newblock {\em J. Sci. Comput.}, 53:102--128, 2012.}

\bibitem{ELW06}
E.~Hairer, C.~Lubich, and G.~Wanner.
\newblock {\em Geometric Numerical Integration: Structure-Preserving Algorithms
  for Ordinary Differential Equations}.
\newblock Springer-Verlag, Berlin, 2nd edition, 2006.

\bibitem{HeND2016}
D.~He.
\newblock Exact solitary solution and a three-level linearly implicit
  conservative finite difference method for the generalized
  {R}osenau-{K}awahara-{RLW} equation with generalized {N}ovikov type
  perturbation.
\newblock {\em Nonlinear Dyn.}, 85:479--498, 2016.

\bibitem{HXH-AMP-2013}
J.~Hu, Y.~Xu, and B.~Hu.
\newblock Conservative linear difference scheme for {Rosenau-KdV} equation.
\newblock {\em Adv. Math. Phys.}, pages 1--7, 2013.

\bibitem{JCWLMaxwell-2017}
C.~Jiang, W.~Cai, Y.~Wang, and H.~Li.
\newblock A novel sixth order energy-conserved method for three-dimensional
  time-domain {M}axwell's equations.
\newblock {\em arXiv:1705.08125}, 2017.

\bibitem{JCCWarXiv2019}
C.~Jiang, J.~Cui, W.~Cai, and Y.~Wang.
\newblock A novel linearized and momentum-preserving fourier pseudo-spectral
  scheme for the {Rosenau-Korteweg de Vries} equation.
\newblock {\em Numer Methods Partial Differential Eq.}, 39:1558--1582, 2023.

\bibitem{JCQSjsc2022}
C.~Jiang, J.~Cui, X.~Qian, and S.~Song.
\newblock High-order linearly implicit structure-preserving exponential
  integrators for the nonlinear {S}chr\"odinger equation.
\newblock {\em J. Sci. Comput.}, 90, 2022,doi.org/10.1007/s10915-021-01739-x.

\bibitem{JWG19}
C.~Jiang, Y.~Wang, and Y.~Gong.
\newblock Arbitrarily high-order energy-preserving schemes for the
  {C}amassa-{H}olm equation.
\newblock {\em Appl. Numer. Math.}, 151:85--97, 2020.

\bibitem{LiCMA2016}
S.~Li.
\newblock Numerical analysis for fourth-order compact conservative difference
  scheme to solve the {3D} {R}osenau-{RLW} equation.
\newblock {\em Comput. Math. Appl.}, 72:2388--2407, 2016.

{\bibitem{LSR-2019-MP}
X.~Li, J.~Shen, and H.~Rui.
\newblock Energy stability and convergence of {SAV} block-centered finite
  difference method for gradient flows.
\newblock {\em Math. Comp.}, 88:2047--2068, 2019.}

\bibitem{OAAKamc08}
K.~Omrani, F.~Abidi, T.~Achouri, and N.~Khiari.
\newblock A new conservative finite difference scheme for the {R}osenau
  equation.
\newblock {\em Appl. Math. Comput.}, 201:35--43, 2008.

\bibitem{PZAMM2012}
X.~Pan and L.~Zhang.
\newblock On the convergence of a conservative numerical scheme for the usual
  {R}osenau-{RLW} equation.
\newblock {\em Appl. Mathe. Model.}, 36:3371--3378, 2012.

\bibitem{QM08}
G.~R.~W. Quispel and D.~I. McLaren.
\newblock A new class of energy-preserving numerical integration methods.
\newblock {\em J. Phys. A: Math. Theor.}, 41:045206, 2008.

\bibitem{RMB-RJP-2014}
P.~Razborova, L.~Moraru, and A.~Biswas.
\newblock Perturbation of dispersive shallow water waves with
  {R}osenau-{KdV-RLW} equation and power law nonlinearity.
\newblock {\em Rom. J. Phys.}, 59:658--676, 2014.

\bibitem{Rosenau1988}
P.~Rosenau.
\newblock Dynamics of dense discrete systems: high order effects.
\newblock {\em Progr. Theor. Phys.}, 79:1028--1042, 1988.

\bibitem{Sanzs88}
J.~M. Sanz-Serna.
\newblock {Runge-Kutta} schemes for {H}amiltonian systems.
\newblock {\em BIT}, 28:877--883, 1988.

\bibitem{SCbook94}
J.~M. Sanz-Serna and M.~Calvo.
\newblock {\em Numerical Hamiltonian Problems}.
\newblock Chapman \& Hall, London, 1994.

\bibitem{ST06}
J.~Shen and T.~Tang.
\newblock {\em Spectral and High-Order Methods with Applications}.
\newblock Science Press, Beijing, 2006.

{\bibitem{ST-SIAM-2018}
J.~Shen and J.~Xu.
\newblock Convergence and error analysis for the scalar auxiliary variable
  ({SAV}) schemes to gradient flows.
\newblock {\em SIAM J. Numer. Anal.}, 56:2895--2912, 2018.}

\bibitem{SXY18}
J.~Shen, J.~Xu, and J.~Yang.
\newblock The scalar auxiliary variable {(SAV)} approach for gradient.
\newblock {\em J. Comput. Phys.}, 353:407--416, 2018.

\bibitem{SXY19siamrev}
J.~Shen, J.~Xu, and J.~Yang.
\newblock A new class of efficient and robust energy stable schemes for
  gradient flows.
\newblock {\em SIAM Rev.}, 61:474--506, 2019.

\bibitem{Tapley-SISC2022}
B.~K. Tapley.
\newblock Geometric integration of {ODEs} using multiple quadratic auxiliary
  variables.
\newblock {\em SIAM J. Sci. Comput.}, 44:A2651--A2668, 2022.

\bibitem{WLWCMA2017}
H.~Wang, S.~Li, and J.~Wang.
\newblock A conservative weighted finite difference scheme for the generalized
  {R}osenau-{RLW} equation.
\newblock {\em Comput. Math. Appl.}, 36:63--78, 2017.

\bibitem{WZJCM2019}
J.~Wang and Q.~Zeng.
\newblock A fourth-order compact and conservative difference scheme for the
  generalized {R}osenau-{K}orteweg de {V}ries equation in two dimensions.
\newblock {\em J. Comput. Math.}, 37:541--555, 2019.

{\bibitem{WHW-JSC-2021}
M.~Wang, Q.~Huang, and C.~Wang.
\newblock A second order accurate scalar auxiliary variable ({SAV}) numerical
  method for the square phase field crystal equation.
\newblock {\em J. Sci. Comput.}, 8, 2021,
  https://doi.org/10.1007/s10915-021-01487-y.}

\bibitem{WDCAM2018}
X.~Wang and W.~Dai.
\newblock A new implicit energy conservative difference scheme with
  fourth-order accuracy for the generalized {R}osenau-{K}awahara-{RLW}
  equation.
\newblock {\em Comput. Appl. Math.}, 37:6560--6581, 2018.

\bibitem{WDYAA2021}
X.~Wang, W.~Dai, and Y.~Yan.
\newblock Numerical analysis of a new conservative scheme for the {2D}
  generalized {R}osenau-{RLW} equation.
\newblock {\em Appl. Anal.}, 100:2564--2580, 2021.

\bibitem{WCCPMMAS2021}
B.~Wongsaijai, P.~Charoensawan, T.~Chaobankoh, and K.~Poochinapan.
\newblock Advance in compact structure-preserving manner to the
  {R}osenau-{K}awahara model of shallow-water wave.
\newblock {\em Math. Meth. Appl. Sci.}, 44:7048--7064, 2021.

\bibitem{YZW17}
X.~Yang, J.~Zhao, and Q.~Wang.
\newblock Numerical approximations for the molecular beam epitaxial growth
  model based on the invariant energy quadratization method.
\newblock {\em J. Comput. Phys.}, 333:104--127, 2017.

{\bibitem{ZWHWY-ANM-2017}
C.~Zhang, H.~Wang, J.~Huang, C.~Wang, and X.~Yue.
\newblock A second order operator splitting numerical scheme for the ``good"
  {B}oussinesq equation.
\newblock {\em Appl. Numer. Math.}, 119:179--193, 2017.}

\bibitem{ZJ2202Arix}
G.~Zhang and C.~Jiang.
\newblock Arbitrary high-order structure-preserving methods for the quantum
  {Z}akharov system.
\newblock {\em arXiv:2202.13052}, 2022.

\bibitem{ZZ-JAM-2013}
M.~Zheng and J.~Zhou.
\newblock An average linear difference scheme for the generalized {Rosenau-KdV}
  equation.
\newblock {\em J. Appl. Math.}, pages 1--9, 2014.

\end{thebibliography}

\end{document}